\documentclass[12pt,twoside]{article}
\usepackage{times}
\usepackage{mathrsfs}
\usepackage{amsmath,amssymb}
\usepackage{color}
\usepackage{fancyhdr}
\allowdisplaybreaks

 \def \nn{\nonumber}

\newcommand{\pf}{\noindent {\bf Proof. \hspace{2mm}}}
\newcommand{\ef}{ \hfill $ \Box $ \vskip 3mm}
\newcommand{\be}{\begin{equation}}
\newcommand{\ee}{\end{equation}}
\newcommand{\bea}{\begin{eqnarray}}
\newcommand{\eea}{\end{eqnarray}}

\newcommand{\bR}{{\mathbb R}}
\newcommand{\bN}{{\mathbb N}}

\def\p{\partial}
\def\la{\lambda}

\def\al{\alpha}

\def\q{\quad}

\def\g{\gamma}

\def\dl{\delta}
\def\ve{\varepsilon}

\def\lt{\left}
\def\rt{\right}

\def\dl{\delta}
\def\Dl{\Delta}
\def\i{\infty}

\def\supp{\text{supp }}

\def\p{\partial}
\def\f{\frac}

\def\al{\alpha}

\def\o{\omega}

\def\q{\quad}
\def\qq{\qquad}
\def\s{\sqrt}
\def\ew{e^{2\psi}}
\def\wh{\hat{w}}
 \allowdisplaybreaks

\begin{document}
\let\oldsection\section
\renewcommand\section{\setcounter{equation}{0}\oldsection}
\renewcommand\thesection{\arabic{section}}
\renewcommand\theequation{\thesection.\arabic{equation}}
\newtheorem{claim}{\noindent Claim}[section]
\newtheorem{theorem}{\noindent Theorem}[section]
\newtheorem{lemma}{\noindent Lemma}[section]
\newtheorem{proposition}{\noindent Proposition}[section]
\newtheorem{definition}{\noindent Definition}[section]
\newtheorem{remark}{\noindent Remark}[section]
\newtheorem{corollary}{\noindent Corollary}[section]
\newtheorem{example}{\noindent Example}[section]

\title{Global existence and asymptotic behavior of solutions to the Euler equations with time-dependent damping}

\author{Xinghong Pan
\footnote{E-mail:math.scrat@gmail.com.}\vspace{0.5cm}\\
\small  (Department of Mathematics and
IMS, Nanjing University, Nanjing 210093, China.)\\
\vspace{0.5cm}
}

\date{}
\maketitle

\centerline {\bf Abstract} \vskip 0.3 true cm

We study the isentropic Euler equations with time-dependent damping, given by $\frac{\mu}{(1+t)^\lambda}\rho u$. Here, $\lambda,\mu$ are two non-negative constants to describe the decay rate of damping with respect to time. We will investigate the global existence and asymptotic behavior of small data solutions to the Euler equations when $0<\lambda<1,0<\mu$ in multi-dimensions $n\geq 1$. The asymptotic behavior will coincide with the one that obtained by many authors in the case $\lambda=0$. We will also show that the solution can only decay polynomially in time while in the three dimensions, the vorticity will decay exponentially fast.

\vskip 0.3 true cm

{\bf Keywords:} Euler equations, time-dependent damping, global existence, asymptotic behavior, vorticity, exponential decay.
\vskip 0.3 true cm

{\bf Mathematical Subject Classification 2010:} 35L70, 35L65, 76N15.

\section{Introduction}
\q This paper deals with the isentropic Euler equations with time-dependent damping in multi-dimensions:
\be
\left\{
\begin{aligned}
&\p_t \rho+\nabla\cdot(\rho u)=0, \\
&\p_t(\rho u)+ \nabla\cdot(\rho u\otimes u)+\nabla p=-\f{\mu}{(1+t)^\la}\rho u,\\
&\rho|_{t=0}=1+\rho_0(x),u|_{t=0}= u_0(x),
\end{aligned}
\right.\label{1.1}
\ee
where $\rho_0(x)\in\bR, u_0(x)\in \bR^n$, supported in $\{x\in\bR^n||x|\leq R\}$. Here $\rho(t,x),u(t,x)$ and $p(t,x)$ represent the density, fluid velocity and pressure respectively and $\la, \mu$ are two positive constants to describe the decay rate of the damping in time. We assume the fluid is a polytropic gas which means we can assume $p(\rho)=\f{1}{\g}\rho^\g, \g>1$.

As is well known, when the damping vanishes, shock will form. For the mathematical analysis of finite-time formation of singularities, readers can see Alinhac\cite{Alinhac}, Chemin\cite{Chemin}, Courant-Friedrichs\cite{CF}, Christodoulou\cite{Christodoulou}, Makino-Ukai-Kawashima\cite{MUK},  Rammaha\cite{Rammaha} as well as Chen-Liu\cite{CL}, Yin\cite{Yin}, Sideris\cite{Sideris1,Sideris2} and references therein for more detail.

The Euler equations with non time-decayed damping are
\be
\left\{
\begin{aligned}
&\p_t \rho+\nabla\cdot(\rho u)=0, \\
&\p_t(\rho u)+ \nabla\cdot(\rho u\otimes u)+\nabla p=-\kappa\rho u,\\
\end{aligned}
\right.\label{1.2}
\ee
where $\kappa$ is the damping constant and $1/\kappa$ can be regarded as the relaxation time of some physical fluid. Many authors have proven the global existence and uniqueness of smooth solutions to system \eqref{1.2} with small data. Also the asymptotic behavior of the smooth solution was studied. For the 1d Euler equations, see Dafermos\cite{Dafermos}, Hsiao-Liu\cite{HL}, Huang-Marcati-Pan\cite{HMP},Nishida\cite{Nishida}, Nishihara-Wang-Yang\cite{NWY} and their references. For the multi-dimensional case, Wang-Yang\cite{WY} give the pointwise estimates of the solution by using some energy methods and estimating the Green function to the linearized system. Sideris-Thomases-Wang\cite{STW} proves a similar result by using a simpler approach. They both proved that the smooth solution decays in maximum norm to the background state at a rate of $(1+t)^{-\f{3}{2}}$ in 3 dimensions.  Sideris-Thomases-Wang also show that the smooth solution has a polynomially decayed lower bound in time while the vorticity will decay exponentially. Tan-Wang\cite{TW}, Jiu-Zheng\cite{JZ} study this problem in the frame of Besov space and obtain similar asymptotic behavior of the solution. Also see Kong-Wang \cite{KW} for extension.

It is natural to ask whether the global solution exists when the damping is decayed and what is the critical decay rate to separate the global existence and the finite-time blow up of solutions with small data. The papers \cite{HY,HWY},\cite{Pan1,Pan2} have done some inspection on this topic where the authors systematically study the case that the damping decays with time as $\f{\mu}{(1+t)^\la}$. They believe that there is a pair of non-negative critical exponent $(\la_c(n),\mu_c(n))$, depending on the space dimension $n$, such that

 {\it when $0\leq\la< \la_c(n),0<\mu$ or $\la=\la_c(n), \mu>\mu_c(n)$, \eqref{1.1} have global existence of small-data solutions; while when $\la= \la_c(n),\mu\leq \mu_c(n)$ and $\la>\la_c(n), 0\leq\mu$, the smooth solutions of \eqref{1.1} will blow up in finite time.}

 In two and three dimensions, Hou-Witt-Yin \cite{HWY}, Hou-Yin \cite{HY} have shown that the critical exponent is $(\la_c(3),\mu_c(3))$$=(1,0)$ and $(\la_c(2),\mu_c(2))$$=(1,1)$; while in one dimension, the critical exponent is $(\la_c(1), \mu_c(1))=(1,2)$, partly presented in \cite{Pan1,Pan2}. 

  This paper deals with the case $0<\la<\la_c(n)=1, 0<\mu$ in multi-dimensions $n\geq 1$. We will obtain the global existence and asymptotic behavior of the solution to system \eqref{1.1}. The method used here will be completely different from that in \cite{HY,HWY}. Also we will show that the convergence rate of the solution to the background constant state $(1,0)$ will coincide with that one obtained by Wang-Yang\cite{WY} and Sideris-Thomases-Wang\cite{STW} when $\la=0$ which indicates that the time-asymptotic behavior of the solution is the diffusion wave of the corresponding linear system.

 The proof of the global existence of the solution is based on the method of weighted energy estimates for symmetric hyperbolic system by introducing the sound speed as a new variable rather than the density. We will establish some weighted a prior estimates to the solution. The choice of the weight is inspired by the corresponding linear wave equation with effective damping satisfied by the sound speed. See \eqref{2.3}. The weight can be found in Nishihara\cite{Nishihara}. The local-existence result stated in Kato\cite{Kato} or Majda\cite{Majda} and the continuity argument can assure the global existence of the solution.

 The estimates of convergence rate of the solution to the background state come from the investigation of the fundamental solutions to the phase function of the corresponding linear wave equation which can be found in Wirth\cite{Wirth}. We will show that $L^2$ and $L^\i$ norms of the solution to system \eqref{1.1} will present a decay estimate similar to that of the corresponding linear dissipative wave equation with the same damping.

 Let $H^l(\bR^n)$ be the usual Sobolev space with its norm
\be
\|f\|_{H^l}\triangleq \sum\limits^l_{k=0}\sum\limits_{|\al|=k}\|\p^\al_xf\|_{L^2}. \nn
\ee
where $\p^\al_x=\p^{\al_1}_1...\p^{\al_n}_n$, $\al=(\al_1,...,\al_n)$. Later for convenience, we will use $\p^k_x= \sum\limits_{|\al|=k}\p^\al_x$
, use $\|\cdot\|_p$ to denote $\|\cdot\|_{L^p}$ and $\|\cdot\|=\|\cdot\|_2$.

We state the global existence result as follows.
\begin{theorem}
Denote $B=\f{(1+\la)n}{2}-\dl$, where $\dl\in(0,\f{(1+\la)n}{2}]$ can be arbitrarily small. Suppose that $n\geq 1$, $0<\la<1, 0<\mu$ and $(\rho_0,u_0)\in H^{s+m}(\bR^n)$, supported in $\{x\in\bR^n||x|\leq R\}$, where $s=[n/2]+1$ and $m\geq2$. Then there exists a $\ve_0=\ve_0(\dl,\la,\mu,R)$ such that for any $0\leq\ve\le\ve_0$, when $\|(\rho_0,u_0)\|_{H^{s+m}}\leq\ve$, there exists a unique global classical solution $(\rho(t,x),u(t,x))$ of \eqref{1.1} satisfying
\bea
&&(1+t)^{B+1+\la}\Big(\|\p_t\rho(t)\|^2_{H^{s+m-1}}+\|\p_x\rho(t)\|^2_{H^{s+m-1}}+\|\p_xu(t)\|^2_{H^{s+m-1}}\Big) \nn  \\
&&+(1+t)^{B+2\la}\Big(\|\p_tu(t)\|^2_{H^{s+m-1}}\Big) \nn  \\
&&+(1+t)^B\Big(\|(\rho(t)-1)\|^2 +\|u(t)\|^2\Big)   \nn  \\
&\leq& C_{\la,\mu,\dl,R}\Big(\|\rho_0\|^2_{H^{s+m}}+\|u_0\|^2_{H^{s+m}}\Big).  \label{1.3}
\eea
\end{theorem}

 Throughout this paper we will denote a generic constant by $C$ which may be different from line to line.
\begin{remark}
From \eqref{1.3}, using the Sobolev embedding, for $k=0,1,...,m-1$, we have
\bea
&&(1+t)^{\f{B+1+\la}{2}}\Big(\|\p^k_x\p_t\rho(t)\|_\i+\|\p^{k+1}_x\rho(t)\|_{\i}+\|\p^{k+1}_xu(t)\|_{\i}\Big) \nn  \\
&&+(1+t)^{\f{B}{2}+\la}\Big(\|\p^k_x\p_tu(t)\|_{\i}\Big) \nn  \\
&&+(1+t)^{\f{B}{2}}\Big(\|(\rho(t)-1)\|_\i +\|u(t)\|_\i\Big)   \nn  \\
&\leq& C\ve.  \label{1.4}
\eea
\end{remark}

From \eqref{1.4}, we see that $\|(\rho(t)-1)\|_\i +\|u(t)\|_\i\leq C\ve(1+t)^{-\f{1+\la}{4}n+\f{\dl}{2}}$. When $\la=0$, it decays slower than what we expect for $(1+t)^{-\f{n}{2}}$ as shown in \cite{WY} and \cite{STW}. So next, based on the investigation of the properties to the corresponding linear wave equation \eqref{3.1} of system \eqref{1.1}, we have the following further asymptotic behavior of the solution.

\begin{theorem}
Define $k_c\triangleq\f{1+\la}{1-\la}(n+1)-n-\f{2\dl}{1-\la}$ and $m\geq k_c+2$. Then under the assumption of Theorem1.1, we have the following asymptotic behavior of the solution $(\rho,u)$ in $L^2$ and $L^\i$ norms.\\
For $\rho$:
\bea
&&\|\p^k_x(\rho-1)\|_\i
\leq C\ve\lt\{
\begin{aligned}
&(1+t)^{-(1-\la)\f{n+k}{2}}\qq\q 0\leq k\leq k_c;\\
&(1+t)^{-(1+\la)\f{n+1}{2}+\dl}\qq k_c\leq k\leq m-2.\\
\end{aligned}
\rt.\label{1.5}\\
&&\|\p^k_x(\rho-1)\|
\leq C\ve\lt\{
\begin{aligned}
&(1+t)^{-(1-\la)(\f{n}{4}+\f{k}{2})}\qq\q 0\leq k\leq k_c+\f{n}{2};\\
&(1+t)^{-(1+\la)\f{n+1}{2}+\dl}\qq\q k_c+\f{n}{2}\leq k\leq s+m-2.\\
\end{aligned}
\rt.\label{1.6}
\eea
While for $u$, due to the damping, it will decay slower than $\rho$ by a factor $(1+t)^\la$. That is
\bea
&&\|\p^k_xu\|_\i
\leq C\ve\lt\{
\begin{aligned}
&(1+t)^{-(1-\la)\f{n+k+1}{2}+\la}\qq\q 0\leq k\leq k_c-1;\\
&(1+t)^{-(1+\la)\f{n+1}{2}+\la+\dl}\qq\q k_c-1\leq k\leq m-3.\\
\end{aligned}
\rt.\label{1.7}\\
&&\|\p^k_xu\|
\leq C\ve\lt\{
\begin{aligned}
&(1+t)^{-(1-\la)(\f{n}{4}+\f{k+1}{2})+\la}\qq\q 0\leq k\leq k_c+\f{n}{2}-1;\\
&(1+t)^{-(1+\la)\f{n+1}{2}+\la+\dl}\qq\q\ \ \ k_c+\f{n}{2}-1\leq k\leq s+m-3.\\
\end{aligned}
\rt.\label{1.8}
\eea
\end{theorem}
\begin{remark}
Noting $\dl>0$ can be arbitrarily small, then $\lim\limits_{\la\rightarrow 0}k_c=1^+$. From \eqref{1.5} and \eqref{1.7}, we have when $\la\rightarrow 0$
\be
\|(\rho-1)(t)\|_\i\leq C(1+t)^{-\f{n}{2}}, \|u(t)\|_\i\leq C(1+t)^{-\f{n+1}{2}}.\nn
\ee
This coincides with the decay rate that obtained in the non-decayed damping case.
\end{remark}

In 3-d Euler equations, denoting $\o=\nabla\times u$, we will derive the exponential decay of $\o$ in time in $L^2$ norm under a positive integration assumption on $\rho_0$. While the solution itself can not decay so fast to its background state. It has a polynomially decayed lower bound. They are presented in the following Theorem.
\begin{theorem}
Suppose $q_0\triangleq\int_{\bR^n}\rho_0dx>0$. Under the assumption of Theorem1.1, there exists a $t_0>0,$ depending on $\la,\mu, R$ such that when $t\geq t_0$,
\be
\|(\rho-1)(t)\|\geq Cq_0(1+t)^{-\f{n}{2}},\qq \|u(t)\|\geq Cq_0(1+t)^{-\f{n+2}{2}}.\label{1.9}
\ee

On the other hand, in three dimensions, the $L^2$ norm of the vorticity $\o$ decays exponentially in time, satisfying the following estimate
\be
\|\o(\cdot,t)\|\leq Ce^{-C(1+t)^{1-\la}}. \label{1.10}
\ee
\end{theorem}
\begin{remark}
The idea of proving Theorem1.3 comes from \cite{STW}, where the authors deal with the case $n=3$ and the damping is non-decayed in time.
\end{remark}

The paper is organized as follows. In Section 2, we reformulate the Euler equations into a symmetric hyperbolic system. Then based on a fundamental weighted energy inequality in Lemma2.1 and using detailed weighted energy estimates, we prove the global existence of smooth solutions with small data. In Section 3, by investigating the structure to the linear system of \eqref{1.1}, we give the asymptotic behavior of the solution. In Section 4, we show that the vorticity decay exponentially while the solution itself has a polynomially decayed lower bound. In the Appendix, we give the proof of the fundamental weighted energy inequality.

\section{Global Existence}

\q In this Section, First we reformulate system \eqref{1.1} to a symmetric system. Remember $c=\s{P'(\rho)}=\rho^{\f{\g-1}{2}}$. First we transform \eqref{1.1} into the following system
\be
\left\{
\begin{aligned}
&\f{2}{\g-1}\p_t c+c\nabla\cdot u+\f{2}{\g-1}u\nabla\cdot c=0, \\
&\p_t u+ u\cdot\nabla u+\f{2}{\g-1}c\nabla c+\f{\mu}{(1+t)^\la} u=0,\\
&c|_{t=0}=1+ c_0(x),u|_{t=0}= u_0(x),
\end{aligned}
\right.\label{2.1}
\ee
where $c_0(x)\in \bR$, supported in $\{x\in\bR^n||x|\leq R\}.$

Let $v=\f{2}{\g-1}(c-1)$, then $(v, u)$ satisfies
\be
\left\{
\begin{aligned}
&\p_t v+\nabla\cdot u=-u\cdot \nabla v-\f{\g-1}{2}v\nabla\cdot u, \\
&\p_t u+ \nabla v+\f{\mu}{(1+t)^\la} u=-u\cdot\nabla u-\f{\g-1}{2}v\nabla v,\\
&v|_{t=0}=v_0(x),u|_{t=0}= u_0(x),
\end{aligned}
\right.\label{2.2}
\ee
where $v_0(x)=\f{2}{\g-1}c_0(x).$
\subsection{A fundamental weighted energy inequality}

From \eqref{2.2}, we have

\be
\p_{tt} v-\Dl v+\f{\mu}{(1+t)^\la}\p_tv=Q(v,u),\label{2.3}
\ee
where
\bea
Q(v,u)&=&\f{\mu}{(1+t)^\la}(-u\cdot\nabla v-\f{\g-1}{2}v\nabla\cdot u)             \nn            \\
&-&\p_t(u\cdot\nabla v-\f{\g-1}{2}v\nabla\cdot u)+\nabla\cdot(u\cdot\nabla u+\f{\g-1}{2}v\nabla v).  \nn
\eea

In the following, we will obtain a fundamental weighted energy inequality about \eqref{2.3}. This technique comes from Nishihara\cite{Nishihara}. Introduce the weight
\be
e^{2\psi},\q \psi(t,x)=a\f{|x|^2}{(1+t)^{1+\la}},  \nn
\ee
where $a=\f{(1+\la)\mu}{8}\Big(1-\f{\dl}{(1+\la)n}\Big)$ and $\dl$ is described in Theorem1.1.
For simplification of notation, we denote
\be
J(t;g)=\int_{\bR^n}\ew g^2(t,x)dx,\q J_\psi(t;g)=\int_{\bR^n}\ew (-\psi_t)g^2(t,x)dx. \nn
\ee
\begin{lemma}
Denote $B=\f{(1+\la)n}{2}-\dl$, where $\dl\in(0,\f{(1+\la)n}{2}]$ can be arbitrarily small. Then the equation \eqref{2.3} has the following weighted energy inequality.
\bea
&&(1+t)^{B+1+\la}\Big[J(t;v_t)+J(t;|\p_xv|)\Big]+(1+t)^{B}J(t;v) \nn  \\
&&\qq+\int^t_0(1+\tau)^{B+1+\la}\Big[J_\psi(\tau;v_\tau)+J_\psi(\tau;|\nabla v|)\Big]d\tau\nn\\
&&\qq+\int^t_0\Big[(1+\tau)^{B+1}J(\tau,v_\tau)+(1+\tau)^{B+\la} J(\tau,|\nabla v|)\Big]d\tau  \nn  \\
&&\qq+\int^t_0\Big[(1+\tau)^{B}J_\psi(\tau,v)+(1+\tau)^{B-1}J(\tau,v)\Big]d\tau \nn  \\
&\leq&C\|(v_0,u_0)\|^2_{H^1}+C G(t), \nn
\eea
where
\be
G(t)\triangleq\int^t_0\int_{\bR^n}\ew \Big\{\big[(1+\tau)^{B+1+\la}v_\tau+ (1+\tau)^{B+\la} v\big] Q(v,u)\Big\}dxd\tau,\nn
\ee
and $C$ depends on $\la,\mu,\dl,R$.
\end{lemma}
\pf We give the proof of Lemma 2.1 in the Appendix. \ef
\subsection{Some a priori weighted energy estimates}

Define the weighted Sobolev space $H^{s+m}_\psi(\bR^n)$ as
\be
H^{s+m}_\psi(\bR^n)=\{f|e^\psi\p^k_xf\in L^2(\bR^n),0\leq k\leq s+m,k\in\bN\}  \nn
\ee
with its norm
\be
\|f\|_{H^{s+m}_\psi}=\sum^{s+m}\limits_{k=0}\Big(\int_{\bR^n}(\p^k_xf)^2\ew dx\Big)^{\f{1}{2}}.\nn
\ee
We define the weighted energy as follows
\bea
E^\psi_{s+m}(T)&=&:\sup\limits_{0<t<T}\lt\{(1+t)^{B+1+\la}\big(\|v_t(t)\|^2_{H^{s+m-1}_\psi}+\|\p_xv(t)\|^2_{H^{s+m-1}_\psi}+\|\p_xu(t)\|^2_{H^{s+m-1}_\psi}\big)\rt.\nn\\
&&\qq\qq+(1+t)^{B}\big(\|v(t)\|_{L^2_\psi}^2+\|u(t)\|_{L^2_\psi}^2\big)\Big\}^{\f{1}{2}}. \nn
\eea
Because the initial data is compactly supported, we have
\bea
E^\psi_{s+m}(0)&\leq& C\Big(\|(v_t(0),\p_xv(0),\p_xu(0))\|_{H^{s+m-1}}+\|(v(0),u(0))\|\Big)\nn\\
&\leq&C\|(v_0,u_0)\|_{H^{s+m}}\leq C\ve. \nn
\eea

In the following, we will estimate $(v,u)$ under the a priori assumption
\be
E^\psi_{s+m}(t)\leq M\ve,  \label{2.4}
\ee
where $M$, independent of $\ve$, will be determined later. By choosing $M$ large and $\ve$ sufficient small, we will prove
\be
E^\psi_{s+m}(t)\leq \f{1}{2}M\ve.  \label{2.5}
\ee

We will first obtain some a prior estimates for the 1-order derivatives. Then the higher derivatives will be handled in a similar way.
From \eqref{2.4}, by Sobolev embedding, we have
\bea
&&\Big\{(1+t)^{\f{B}{2}}\|(v(t),u(t))\|_{L^\i}+(1+t)^{\f{B+1+\la}{2}}\|(v_t(t),\p_xv(t),\p_xu(t))\|_{L^\i}\Big\}\nn\\
&\leq& CE^\psi_{[n/2]+2}(t)\leq CM\ve.  \label{2.6}
\eea
From $\eqref{2.2}_2$, we see
\bea
&&(1+t)^{B/2+\la}\|u_t(t)\|_{L^\i}\nn\\
&\leq&  (1+t)^{B/2+\la}\Big\{\|\nabla v(t)\|_{L^\i}+(1+t)^{-\la}\|u (t)\|_{L^\i}+\|u\cdot\nabla u\|_{L^\i}+\|v\nabla v\|_{L^\i}\Big\}\nn\\
&\leq&CM\ve.  \label{2.7}
\eea
Now for simplification of notation in our later computation, we introduce
\bea
&&P_0(t)\nn\\
&\triangleq&P(t;v,u,v_t,u_t,\p_x v,\p_x u)\nn\\
&\triangleq&(1+t)^{B+1+\la}\Big[J_\psi(t;v_t)+J_\psi(t;|\p_x v|)+J_\psi(t;|\p_x u|)\Big]\nn\\
&&+(1+t)^{B+1}\big[J(t;v_t)+J(t;|\p_x u|)\big]+(1+t)^{B+\la}\big[J(t;|u_t|)+J(t;|\p_x v|)\big]\nn\\
&&+(1+t)^{B}\big[J_\psi(t;v)+J_\psi(t;| u|)\big]\nn\\
&&+(1+t)^{B-\la}J(t;|u|)+(1+t)^{B-1}J(t; v).\label{2.8}
\eea

In the following, we will give three estimates for the first order derivatives under the a prior assumption, which can reflect our idea of proving the   \eqref{2.5}.\\

\noindent \textbf{Estimate 1}

\begin{lemma}
Under the a prior assumption \eqref{2.4}, we have the following estimate
\bea
G(t)&\leq&  CM\ve(1+t)^{B+1+\la}\big[J(t;|v_t|)+J(t;|\p_x v|)\big]\nn\\
      &&+C\ve \Big(E^\psi_1(0)\Big)^2+CM\ve\int^t_0P_0(\tau)d\tau,\label{2.9}
\eea
where $P_0(\tau)$ is defined in \eqref{2.8}.
\end{lemma}
\pf
Remember the formulation of $Q(v,u)$, we have
\bea
G(t)
&=&\int^t_0\int_{\bR^n}\ew \Big((1+\tau)^{B+1+\la}v_\tau+ (1+\tau)^{B+\la} v\Big)\times\nn\\
&&\qq\qq\Big\{\underbrace{\f{\mu}{(1+\tau)^\la}\Big(-u\cdot\nabla v-\f{\g-1}{2}v\nabla\cdot u\Big)}_{I_1}\nn \\
&&\qq\qq\q\underbrace{-\Big(\p_\tau(u\cdot\nabla v)+\f{\g-1}{2}\p_\tau(v\nabla\cdot u)\Big)}_{I_2} \nn \\
&&\qq\qq\q\underbrace{+\Big(\nabla\cdot(u\cdot\nabla u)+\f{\g-1}{2}\nabla\cdot(v\nabla v)\Big)}_{I_3}\Big\}dxd\tau.  \nn
\eea
Now we estimate $I_i(i=1,2,3)$ term by term under the a priori estimate \eqref{2.4}.\\
Using \eqref{2.6}, \eqref{2.7} and Cauchy-Schwartz inequality, we have
\bea
I_1&=&\mu\int^t_0\int_{\bR^n}\ew(1+\tau)^{B+1}v_\tau\Big(-u\cdot\nabla v-\f{\g-1}{2}v\nabla\cdot u\Big)dxd\tau  \nn\\
 &&+\mu\int^t_0\int_{\bR^n}(1+\tau)^{B}\ew v\Big(-u\cdot\nabla v-\f{\g-1}{2}v\nabla\cdot u\Big)dxd\tau \nn  \\
 &\leq&C\|(1+\tau)^{\f{1+\la}{2}}\nabla v(\tau)\|_\i\int^t_0\Big[(1+\tau)^{B+1}J(\tau;v_\tau)+(1+\tau)^{B-\la}J(\tau;|u|)\Big]d\tau \nn  \\
 &&+C\| v(\tau)\|_\i\int^t_0\Big[(1+\tau)^{B+1}J(\tau;v_\tau)+(1+\tau)^{B+1}J(\tau;|\nabla u|)\Big]d\tau \nn  \\
 &&+C\|v(\tau)\|_\i\int^t_0\Big[(1+\tau)^{B-\la}J(\tau;|u|)+(1+\tau)^{B+\la}J(\tau;|\nabla v|)\Big]d\tau \nn  \\
 &&+C\|v(\tau)\|_\i\int^t_0\Big[(1+\tau)^{B-1}J(\tau;v)+(1+\tau)^{B-1}J(\tau;|\nabla u|)\Big]d\tau \nn  \\
 &\leq&CM\ve \int^t_0P_0(\tau)d\tau.  \label{2.10}
\eea
And
\bea
I_2&=&-\int^t_0\int_{\bR^n}\ew \Big((1+\tau)^{B+1+\la}v_\tau+ (1+\tau)^{B+\la} v\Big)\times\nn\\
&&\qq\qq\qq\Big\{\big(u_\tau\cdot\nabla v+\f{\g-1}{2}v_\tau\nabla\cdot u\big) \nn \\
&&\qq\qq\qq\q+\underbrace{\Big(u\cdot\nabla v_\tau\Big)}_{I_{2,1}}+\underbrace{\Big(\f{\g-1}{2}v\nabla\cdot u_\tau\Big)}_{I_{2,2}}\Big\}dxd\tau\nn\\
&\leq&CM\ve\int^t_0P_0(\tau)d\tau+I_{2,1}+I_{2,2}.    \label{2.11}
\eea
Next we use integration by parts to estimate $I_{2,1}$ and $I_{2,2}$.
\bea
I_{2,1}&=&-\int^t_0\int_{\bR^n}\ew \Big((1+\tau)^{B+1+\la}v_\tau+ (1+\tau)^{B+\la} v\Big)\Big(u\cdot\nabla v_\tau\Big)dxd\tau\nn\\
&=-&\f{1}{2}\int^t_0\int_{\bR^n}\ew (1+\tau)^{B+1+\la}u\cdot \nabla v^2_\tau dxd\tau-\int^t_0\int_{\bR^n}\ew(1+\tau)^{B+\la} vu\cdot\nabla v_\tau dxd\tau\nn\\
&=&\f{1}{2}\int^t_0\int_{\bR^n}\ew (1+\tau)^{B+1+\la}v^2_\tau(\nabla\cdot u+\underbrace{2u\cdot\nabla\psi}_{I^1_{2,1}})dxd\tau  \nn \\
&&+\int^t_0\int_{\bR^n}\ew(1+\tau)^{B+\la} v_\tau(v\nabla\cdot u+u\cdot\nabla v+\underbrace{2vu\cdot\nabla\psi}_{I^2_{2,1}})dxd\tau \nn  \\
&\leq&CM\ve\int^t_0P_0(\tau)d\tau+I^1_{2,1}+ I^2_{2,1}.    \label{2.12}
\eea
Noting \eqref{5.2}, we have
\be
|\nabla\psi|\leq C(1+t)^{-\f{\la}{2}}(-\psi_t)^{\f{1}{2}}.  \label{2.13}
\ee
Using Cauchy-Schwartz inequality, we obtain
\bea
&&|I^1_{2,1}|+ |I^2_{2,1}|\nn\\
&\leq& C\|(1+\tau)^\f{1+\la}{2}v_\tau\|_{L^\i}\int^t_0\big[(1+\tau)^{B+1}J(\tau;v_\tau)+(1+\tau)^BJ_\psi(\tau;|u|)\big]dxd\tau\nn\\
&\leq&C\|(1+\tau)^\f{1+\la}{2}v_\tau\|_{L^\i}\int^t_0\big[(1+\tau)^{B-1}J(\tau;v)+(1+\tau)^BJ_\psi(\tau;|u|)\big]dxd\tau\nn\\
&\leq&CM\ve\int^t_0P_0(\tau)d\tau. \label{2.14}
\eea
Combing \eqref{2.12} and \eqref{2.14}, we have
\be
|I_{2,1}|\leq CM\ve\int^t_0P_0(\tau)d\tau.\label{2.15}
\ee
For $I_{2,2}$, we have
\bea
I_{2,2}&=&\int^t_0\int_{\bR^n}\ew \Big((1+\tau)^{B+1+\la}v_\tau+ (1+\tau)^{B+\la} v\Big)\Big(\f{\g-1}{2}v\nabla\cdot u_\tau\Big)dxd\tau \nn \\
       &=&C\underbrace{\int^t_0\int_{\bR^n}\ew (1+\tau)^{B+\la}v^2\nabla\cdot u_\tau dxd\tau}_{I^1_{2,2}}  \nn  \\
       &&+C\underbrace{\int^t_0\int_{\bR^n}\ew (1+\tau)^{B+1+\la}v_\tau v\nabla\cdot u_\tau dxd\tau}_{I^2_{2,2}}. \label{2.16}
\eea
Using integration by parts, \eqref{2.13} and Cauchy-Schwartz inequality, we have
\bea
I^1_{2,2}&=&-C\int^t_0\int_{\bR^n}\ew (1+\tau)^{B+\la}u_\tau\cdot (2v\nabla\cdot v+2v^2\nabla\psi)  dxd\tau \nn \\
         &\leq&CM\ve\int^t_0P_0(\tau)d\tau.\label{2.17}
\eea
For the estimate of $I^2_{2,2}$, we use \eqref{2.2} and integration by parts.
From $\eqref{2.2}_1$, we have
\be
\nabla\cdot u=-\f{v_t+u\cdot\nabla v}{1+\f{\g-1}{2}v}, \label{2.18}
\ee
\be
\nabla\cdot u_t=-\f{v_{tt}+u_t\cdot\nabla v+u\cdot\nabla v_t}{1+\f{\g-1}{2}v}+\f{\f{\g-1}{2}v_t(v_t+u\cdot\nabla v)}{(1+\f{\g-1}{2}v)^2}. \label{2.19}
\ee
From \eqref{2.6}, we have
\be
\f{1}{1+\f{\g-1}{2}v}\leq \f{1}{1-CM\ve}\leq 2. \label{2.20}
\ee
Inserting \eqref{2.19} and \eqref{2.20} into $I^2_{2,2}$, we have
\bea
I^2_{2,2}&=&\int^t_0\int_{\bR^n}\ew (1+\tau)^{B+1+\la}v_\tau v\nabla\cdot u_\tau dxd\tau\nn\\
&=&\int^t_0\int_{\bR^n}\ew (1+\tau)^{B+1+\la}v_\tau v\Big\{-\f{v_{\tau\tau}+u_\tau\cdot\nabla v+u\cdot\nabla v_\tau}
{1+\f{\g-1}{2}v}\nn\\
&&\qq\qq\qq\qq\qq\qq\qq+\f{\f{\g-1}{2}v_\tau(v_\tau+u\cdot\nabla v)}{(1+\f{\g-1}{2}v)^2}\Big\}dxd\tau\nn\\
&\leq&\underbrace{-\int^t_0\int_{\bR^n}\ew(1+\tau)^{B+1+\la}\f{vv_\tau v_{\tau\tau}}{1+\f{\g-1}{2}v}dxd\tau}_{I^{2,1}_{2,2}}\nn\\
&&\underbrace{-\int^t_0\int_{\bR^n}\ew(1+\tau)^{B+1+\la}\f{vv_\tau u\cdot\nabla v_\tau}{1+\f{\g-1}{2}v} dxd\tau}_{I^{2,2}_{2,2}}\nn\\
&&+CM\ve\int^t_0P_0(\tau)d\tau. \label{2.21}
\eea
Using integration by parts in time and \eqref{2.20}, we have
\bea
I^{2,1}_{2,2}&=&-\int^t_0\int_{\bR^n}\ew(1+\tau)^{B+1+\la}\f{v (v^2_{\tau})_\tau}{2\big(1+\f{\g-1}{2}v\big)}dxd\tau\nn\\
&=&-\int_{\bR^n}\ew(1+\tau)^{B+1+\la}\f{v v^2_{\tau}}{2\big(1+\f{\g-1}{2}v\big)}\bigg|^{\tau=t}_{\tau=0}dx\nn\\
&&+\int^t_0\int_{\bR^n}\ew v^2_\tau\bigg\{(1+\tau)^{B+1+\la}\p_\tau\bigg(\f{v}{2\big(1+\f{\g-1}{2}v\big)}\bigg)\nn\\
&&\qq\qq\qq\qq+(B+1+\la)(1+\tau)^{B+\la}\f{v}{2\big(1+\f{\g-1}{2}v\big)}\nn\\
&&\qq\qq\qq\qq+2\psi_t(1+\tau)^{B+1+\la}\f{v}{2\big(1+\f{\g-1}{2}v\big)}\bigg\}dxd\tau\nn\\
&\leq& CM\ve(1+t)^{B+1+\la}J(t;v_t)+C\ve \Big(E^\psi_1(0)\Big)^2+CM\ve\int^t_0P_0(\tau)d\tau.  \label{2.22}
\eea
And
\bea
I^{2,2}_{2,2}&=&-\int^t_0\int_{\bR^n}\ew(1+\tau)^{B+1+\la}\f{v}{2\big(1+\f{\g-1}{2}v\big)}u\cdot\nabla v^2_\tau dxd\tau\nn\\
&=&\int^t_0\int_{\bR^n}\ew(1+\tau)^{B+1+\la}v^2_\tau\bigg\{u\cdot\nabla\Big(\f{v}{2\big(1+\f{\g-1}{2}v\big)}\Big)\nn\\
&&\qq\qq\qq\qq\qq+\f{v}{2\big(1+\f{\g-1}{2}v\big)}\nabla\cdot u\nn\\
&&\qq\qq\qq\qq\qq+\f{v}{1+\f{\g-1}{2}v}u\nabla\psi\bigg\}dxd\tau\nn\\
&\leq&CM\ve\int^t_0P_0(\tau)d\tau. \label{2.23}
\eea
From \eqref{2.21}, \eqref{2.22} and \eqref{2.23}, we have
\be
I^2_{2,2}\leq CM\ve(1+t)^{B+1+\la}J(t;v_t)+C\ve \Big(E^\psi_1(0)\Big)^2+CM\ve\int^t_0P_0\tau)d\tau. \nn
\ee
Summing all the estimates about $I_2$, we can get
\be
I_2\leq CM\ve(1+t)^{B+1+\la}J(t;v_t)+C\ve \Big(E^\psi_1(0)\Big)^2+CM\ve\int^t_0P_0(\tau)d\tau. \label{2.24}
\ee
The estimate of $I_3$ will be essentially the same with $I_2$, actually we can get
\be
I_3\leq CM\ve(1+t)^{B+1+\la}J(t;|\p_xv|)+C\ve \Big(E^\psi_1(0)\Big)^2+CM\ve\int^t_0P_0(\tau)d\tau.\label{2.25}
\ee
Combing the estiamtes \eqref{2.10}, \eqref{2.24} and \eqref{2.25}, we get the estimate
\bea
G(t)&\leq&  CM\ve(1+t)^{B+1+\la}\big[J(t;|v_t|)+J(t;|\p_x v|)\big]\nn\\
      &&+C\ve\Big(E^\psi_1(0)\Big)^2+CM\ve\int^t_0P_0(\tau)d\tau,\nn
\eea
which finishes the proof of the Lemma. \ef
Combing Lemma2.1 and Lemma2.2, we have
\bea
&&(1+t)^{B+1+\la}\Big[J(t;v_t)+J(t;|\nabla v|)\Big]+(1+t)^{B}J(t;v) \nn  \\
&&\qq+\int^t_0(1+\tau)^{B+1+\la}\Big[J_\psi(\tau;v_\tau)+J_\psi(\tau;|\nabla v|)\Big]d\tau\nn\\
&&\qq+\int^t_0\Big[(1+\tau)^{B+1}J(\tau,v_\tau)+(1+\tau)^{B+\la} J(\tau,|\nabla v|)\Big]d\tau  \nn  \\
&&\qq+\int^t_0\Big[(1+\tau)^{B}J_\psi(\tau,v)+(1+\tau)^{B-1}J(\tau,v)\Big]d\tau \nn  \\
&\leq&C \Big(E^\psi_1(0)\Big)^2+CM\ve\int^t_0P_0(\tau)d\tau. \label{2.26}
\eea
\noindent \textbf{Estimate 2}\\

 Next we consider the weighted $L^2$ norm of $(v,u)$. We have the following Lemma.
\begin{lemma}
Under the a prior assumption \eqref{2.4}, we have the following estimate
\bea
&&(1+t)^B\big[J(t;v)+J(t;|u|)]\nn\\
&&+\int^t_0(1+\tau)^BJ_\psi(\tau;|u|)d\tau+\int^t_0(1+\tau)^{B-\la}J(\tau;|u|)d\tau\nn\\
&&-C\int^t_0\big[(1+\tau)^BJ_\psi(\tau; v)+(1+\tau)^{B-1}J(\tau;v)\big]d\tau\nn\\
&\leq& C\big(J(0;v)+J(0;|u|)\big)+CM\ve\int^t_0P_0(\tau)d\tau.\label{2.27}
\eea
where $P_0(\tau)$ is defined in \eqref{2.8}.
\end{lemma}
\pf Multiplying $\eqref{2.2}_2$ by $(K+t)^Be^{2\psi}u$ yields
\bea
&&\p_t\Big[(K+t)^B\f{\ew}{2}|u|^2\Big]+\ew(K+t)^B(-\psi_t)|u|^2\nn\\
&&\Big(\f{\mu(K+t)^B}{(1+t)^\la}-\f{B(K+t)^{B-1}}{2}\Big)\ew |u|^2+(K+t)^B\ew u\cdot\nabla v\nn\\
&=&(K+t)^B\ew u\cdot(-u\cdot\nabla u-\f{\g-1}{2}v\nabla v).   \label{2.28}
\eea
Integrating \eqref{2.28} on $\bR^n\times [0,t]$ and choosing $K$ be large, we  can get
\bea
&&(K+t)^BJ(t;|u|)+\int^t_0(K+\tau)^BJ_\psi(\tau;|u|)d\tau+\int^t_0(K+\tau)^{B-\la}J(\tau;|u|)d\tau\nn\\
&&+C\underbrace{\int^t_0\int_{\bR^n}\ew (K+\tau)^Bu\cdot\nabla vdxd\tau}_{I_1}\nn\\
&\leq&C\underbrace{\int^t_0\int_{\bR^n}(K+\tau)^B\ew u\cdot(-u\cdot\nabla u-\f{\g-1}{2}v\nabla v)dxd\tau}_{I_2}+CJ(0;|u|). \label{2.29}
\eea
We come to estimate $I_1$ and $I_2$.\\
Using integration by parts and $\eqref{2.2}_1$, we have
\bea
I_1&=&-\int^t_0\int_{\bR^n}(K+\tau)^B\ew v(\nabla\cdot u+2u\cdot\nabla\psi)dxd\tau\nn\\
&=&\int^t_0\int_{\bR^n}(K+\tau)^B\ew v (v_\tau+u\cdot\nabla v+\f{\g-1}{2}v\nabla\cdot u-2u\cdot\nabla\psi)dxd\tau\nn\\
&=&\underbrace{\int^t_0\int_{\bR^n}\Big\{\p_\tau\Big[(K+\tau)^B\ew\f{v^2}{2}\Big]+(K+\tau)^B\ew(-\psi_\tau)v^2-B/2(K+\tau)^{B-1}\ew v^2\Big\}dxd\tau}_{I_{1,1}}\nn\\
&&+\underbrace{\int^t_0\int_{\bR^n}(K+\tau)^B\ew v(u\cdot\nabla v+\f{\g-1}{2}v\nabla\cdot u)dxd\tau}_{I_{1,2}}\nn\\
&&-2\underbrace{\int^t_0\int_{\bR^n}(K+\tau)^B\ew vu\cdot\nabla\psi dxd\tau}_{I_{1,3}}. \label{2.30}
\eea
We see that
\bea
I_{1,1}&=&1/2(K+t)^BJ(t;v)-CJ(0;v)\nn\\
&&+\int^t_0(K+\tau)^BJ_\psi(\tau;v)d\tau-B/2\int^t_0(K+\tau)^{B-1}J(\tau;v)d\tau. \label{2.31}
\eea
Using \eqref{2.6} and \eqref{2.7}, we have
\bea
|I_{1,2}|&\leq&\|(K+t)^{(1+\la)/2}\p_x v(t)\|_{L^\i}\int^t_0\big[(K+t)^{B-\la}J(\tau;u)+(K+\tau)^{B-1}J(\tau;v)\big]d\tau\nn\\
&&+\|(K+t)^{(B+1)/2}\p_x u(t)\|_{L^\i}\int^t_0(K+t)^{B-1}J(\tau;v)d\tau\nn\\
&\leq&CM\ve\int^t_0P_0(\tau)d\tau.\label{2.32}
\eea
Using Cauchy-Schwartz inequality and \eqref{5.2}, we obtain
\bea
|I_{1,3}|&\leq&\nu \int^t_0(K+t)^{B-\la}J(\tau;|u|)d\tau+C_\nu\int^t_0(K+t)^BJ_\psi(\tau; v)d\tau,\label{2.33}
\eea
where $\nu$ is a small number.
Combing $\eqref{2.30}-\eqref{2.33}$, we have
\bea
I_{1}&\geq&1/2(K+t)^BJ(t;v)-CJ(0;v)\nn\\
&&+\int^t_0(K+\tau)^BJ_\psi(\tau;v)d\tau-B/2\int^t_0(K+\tau)^{B-1}J(\tau;v)d\tau\nn\\
&&-\nu \int^t_0(K+t)^{B-\la}J(\tau;|u|)d\tau-C_\nu\int^t_0(K+t)^BJ_\psi(\tau; v)d\tau\nn\\
&&-CM\ve\int^t_0P_0(\tau)d\tau.  \label{2.34}
\eea
$I_2$ can be done the same with $I_{1,2}$, so, we have
\be
|I_2|\leq CM\ve\int^t_0P_0(\tau)d\tau.\label{2.35}
\ee
Combing \eqref{2.29}, \eqref{2.34} and \eqref{2.35}, by choosing $\nu$ small, we have
\bea
&&(K+t)^B\big[J(t;v)+J(t;|u|)]\nn\\
&&+\int^t_0(K+\tau)^BJ_\psi(\tau;|u|)d\tau+\int^t_0(K+\tau)^{B-\la}J(\tau;|u|)d\tau\nn\\
&&-C\int^t_0\big[(K+t)^BJ_\psi(\tau; v)+(K+t)^{B-1}J(t;v)\big]d\tau\nn\\
&\leq& C\big(J(0;v)+J(0;|u|)\big)+CM\ve\int^t_0P_0(\tau)d\tau.\nn
\eea
This finishes the proof of Lemma2.3.\ef
\noindent \textbf{Estimate 3}\\

 Next we consider the weighted $L^2$ norm of $(\p_xv,\p_xu)$. We have the following lemma.
\begin{lemma}
Under the a prior assumption \eqref{2.4}, we have the following estimate
\bea
&&(1+t)^{B+1+\la}\big[J(t;|\p_xv|)+J(t;|\p_xu|)]\nn\\
&&+\int^t_0(1+\tau)^{B+1+\la}J_\psi(\tau;|\p_xu|)d\tau+\int^t_0(1+\tau)^{B+1}J(\tau;|\p_xu|)d\tau\nn\\
&&-C\int^t_0\big[(1+\tau)^{B+1+\la}J_\psi(\tau; |\p_xv|)+(1+\tau)^{B+\la}J(\tau;|\p_x v|)\big]d\tau\nn\\
&\leq& C\Big(E^\psi_1(0)\Big)^2+CM\ve\int^t_0P_0(\tau)d\tau,\label{2.36}
\eea
where $P_0(\tau)$ is defined in \eqref{2.8}.
\end{lemma}

\pf
By differentiating $\eqref{2.2}_2$ with respect to $x$ and then multiplying it with $(K+t)^{B+1+\la}\ew\p_xu$, we get
\bea
&&\p_t\Big[(K+t)^{B+1+\la}\ew\f{|\p_xu|^2}{2}\Big]+\ew(K+t)^{B+1+\la}(-\psi_t)|\p_x u|^2\nn\\
&&\q+\bigg[\Big(\mu\f{(K+t)^\la}{(1+t)^\la}-\f{B+1+\la}{2}(K+t)^{\la-1}\Big)(K+t)^{B+1}\ew|\p_x u|^2\bigg]\nn\\
&&\q+(K+t)^{B+1+\la}\ew \p_x u\cdot\nabla\p_xv\nn\\
&=&(K+t)^{B+1+\la}\ew\p_x u\cdot\p_x(-u\cdot\nabla u-\f{\g-1}{2}v\nabla v). \label{2.37}
\eea

By choosing $K$ large and integrating \eqref{2.37} on $\bR^n\times[0,t]$, we can obtain the estimate \eqref{2.36} by using the same trick as that in the proof of Lemma2.3. We omit the detail. \ef
\noindent \textbf{Estimate 4}\\
\\
From $\eqref{2.2}_2$, we can easily have
\bea
&&\int^t_0(1+\tau)^{B+\la}J(\tau;|u_\tau|)d\tau\nn\\
&\leq&\int^t_0(1+\tau)^{B+\la}J(\tau;|\p_x v|)d\tau+\int^t_0(1+\tau)^{B-\la}J(\tau;|u|)d\tau+CM\ve\int^t_0P_0(\tau)d\tau.\nn\\
&&\label{2.38}
\eea
Adding \eqref{2.26} to $\nu_1\cdot\eqref{2.27}+\nu_2\cdot\eqref{2.36}+\nu_3\cdot\eqref{2.38}$, where $\nu_1,\nu_2\ll1$ and $\nu_3\ll\nu_1$, we can get
\bea
&&\Big(E^\psi_1(t)\Big)^2+\int^t_0P_0(\tau)d\tau\nn\\
&\leq&\Big(E^\psi_1(0)\Big)^2+CM\ve\int^t_0P_0(\tau)d\tau.\label{2.39}
\eea
\subsection{Proof of Theorem1.1}

\pf Actually the estimate \eqref{2.39} can be obtained for higher derivatives. Define $P_k(t)$ the same as $P_0$.
\bea
&&P_k(t)\nn\\
&\triangleq&P(t;\p^k_xv,\p^k_x u,\p^k_xv_t,\p^k_xu_t,\p_x \p^k_xv,\p_x \p^k_xu)\nn\\
&\triangleq&(1+t)^{B+1+\la}\Big[J_\psi(t;\p^k_xv_t)+J_\psi(t;|\p_x \p^k_xv|)+J_\psi(t;|\p_x \p^k_xu|)\Big]\nn\\
&&+(1+t)^{B+1}\big[J(t;\p^k_xv_t)+J(t;|\p_x \p^k_xu|)\big]+(1+t)^{B+\la}\big[J(t;|\p^k_xu_t|)+J(t;|\p_x \p^k_xv|)\big]\nn\\
&&+(1+t)^{B}\big[J_\psi(t;\p^k_xv)+J_\psi(t;| \p^k_xu|)\big]\nn\\
&&+(1+t)^{B-\la}J(t;|\p^k_xu|)+(1+t)^{B-1}J(t; \p^k_xv).\nn
\eea
For any $l\in\bN$, differentiating \eqref{2.3} $l$ times with respect to space variable $x$, we have
\be
\p_{tt} \p^l_xv-\Dl \p^l_xv+\f{\mu}{(1+t)^\la}\p_t\p^l_xv=\p^l_xQ(v,u).\label{2.40}
\ee

Then do the same estimate as shown in \textbf{Estimate 1}. We can get the following similar estimate to \eqref{2.26}
\bea
&&(1+t)^{B+1+\la}\big[J(t;\p^l_xv_t)+J(t;|\p_x\p^l_xv|)\big] \nn  \\
&&+\int^t_0(1+\tau)^{B+1+\la}\big[J_\psi(\tau;\p^l_xv_\tau)+J_\psi(\tau;(|\p_x\p^l_xv|)\big]d\tau\nn\\
&&+\int^t_0(1+\tau)^{B+1}J(\tau;\p^l_xv_\tau)+(1+\tau)^{B+\la} J(\tau;|\p_x \p^l_xv|)d\tau  \nn  \\
&&+\int^t_0(1+\tau)^BJ_\psi(\tau,\p^l_xv)+(1+\tau)^{B-1}J(\tau;\p^l_xv)d\tau \nn  \\
&\leq&CE^\psi_{l+1}(0)+CM\ve\sum^l\limits_{k=0}\int^t_0P_k(\tau)d\tau. \label{2.41}
\eea

Differentiating $\eqref{2.2}_2$ $l$ times  with respect to $x$, multiplying it by $(K+t)^{B}\ew\p^l_x u$ and doing the same estimate as in \textbf{Estimate 2}, we have
\bea
&&(1+t)^B\big[J(t;\p^l_xv)+J(t;|\p^l_xu|)]\nn\\
&&+\int^t_0(1+\tau)^BJ_\psi(\tau;|\p^l_xu|)d\tau+\int^t_0(1+\tau)^{B-\la}J(\tau;|\p^l_xu|)d\tau\nn\\
&&-C\int^t_0\big[(1+\tau)^BJ_\psi(\tau; \p^l_xv)+(1+\tau)^{B-1}J(\tau;\p^l_xv)\big]d\tau\nn\\
&\leq&CE^\psi_l(0)+ CM\ve\sum^{l-1}\limits_{k=0}\int^t_0P_k(\tau)d\tau.\label{2.42}
\eea

Differentiating $\eqref{2.2}_2$ $l+1$ times  with respect to $x$, multiplying it by $(K+t)^{B+1+\la}\ew\p^{l+1}_x u$ and doing the same estimate as in \textbf{Estimate 3}, we have
\bea
&&(1+t)^{B+1+\la}\big[J(t;|\p_x\p^l_xv|)+J(t;|\p_x \p^l_xu|)]\nn\\
&&+\int^t_0(1+\tau)^{B+1+\la}J_\psi(\tau;|\p_x \p^l_xu|)d\tau+\int^t_0(1+\tau)^{B+1}J(\tau;|\p_x \p^l_xu|)d\tau\nn\\
&&-C\int^t_0\big[(1+\tau)^{B+1+\la}J_\psi(\tau; |\p_x\p^l_xv|)+(1+\tau)^{B+\la}J(\tau;|\p_x \p^l_xv|)\big]d\tau\nn\\
&\leq&CE^\psi_{l+1}(0)+ CM\ve\sum^l\limits_{k=0}\int^t_0P_k(\tau)d\tau.\label{2.43}
\eea
From $\eqref{2.2}_2$, we also have
\bea
&&\int^t_0(1+\tau)^{B+\la}J(\tau;|\p^l_xu_\tau|)d\tau\nn\\
&\leq&\int^t_0(1+\tau)^{B+\la}J(\tau;|\p_x \p^l_xv|)d\tau+\int^t_0(1+\tau)^{B-\la}J(\tau;|\p^l_xu|)d\tau\nn\\
&&+CM\ve\sum^l\limits_{k=0}\int^t_0P_k(\tau)d\tau.\label{2.44}
\eea
Adding $\eqref{2.41}-\eqref{2.44}$ together the same as \eqref{2.39} for $l$ from $0$ to $s+m-1$ together. Then there exists a $C_0$, depending on $\la,\mu, \dl$ and $R$ such that
\bea
&&\Big(E^\psi_{s+m}(t)\Big)^2+\sum^{s+m-1}\limits_{k=0}\int^t_0P_k(\tau)d\tau\nn\\
&\leq& C\Big(E^\psi_{s+m}(0)\Big)^2+CM\ve\sum^{s+m-1}\limits_{k=0}\int^t_0P_k(\tau)d\tau. \nn\\
&\leq& C^2_0\ve^2+C_0M\ve\sum^{s+m-1}\limits_{k=0}\int^t_0P_k(\tau)d\tau. \nn
\eea
Let $M=2C_0$ and $C_0M\ve<1$. Then we obtain
\bea
E^\psi_{s+m}(t)\leq   \f{1}{2}M\ve.  \nn
\eea

The local existence of symmetrizable hyperbolic equations have been proved by using the fixed point theorem. In order to get the global existence of the system, we only need a priori estimate. Based on the a priori estimate \eqref{2.4}, we yield \eqref{2.5}. This finishes the proof of Theorem 1.1.\ef

\section{Asymptotic Behavior of the Solution}

\q In this Section, we restudy the decay rate of the smooth solution to system \eqref{1.1} and obtain the asymptotic behavior of the solution to the background state $(1,0)^{\tau}$. For this purpose, we consider the linear equation of \eqref{2.3}
\be
w_{tt}-\Dl w+\f{\mu}{(1+t)^\la}w_t=f. \label{3.1}
\ee

This equation involves in the study of the Cauchy problem to the wave equation with time-dependent damping.
\be
\lt\{
\begin{aligned}
&w_{tt}-\Dl w+\f{\mu}{(1+t)^\la}w_t=f,\\
&v(0,x)=w_0(x),\q w_t(0,x)=w_1(x).
\end{aligned}
\rt.\label{3.2}
\ee
where $0<\la<1, 0<\mu$.
\subsection{Decay rate of solutions to the linear wave equation}

\q If we apply the partial Fourier transform with respect to the space variables to \eqref{3.2}, we can get
\be
\lt\{
\begin{aligned}
&\wh_{tt}+|\xi|^2 \wh+\f{\mu}{(1+t)^\la}\wh_t=\hat{f},\\
&\wh(0,\xi)=\wh_0(\xi),\q \wh_t(0,\xi)=\wh_1(\xi),
\end{aligned}
\rt.\label{3.3}
\ee
where $\xi$ is the frequency parameter. The solution of \eqref{3.3} can be represented by the sum of the solution of the following two equations
\be
\lt\{
\begin{aligned}
&\wh^1_{tt}+|\xi|^2 \wh^1+\f{\mu}{(1+t)^\la}\wh^1_t=0,\\
&\wh(0,\xi)=\wh_0(\xi),\q \wh_t(0,\xi)=\wh_1(\xi),
\end{aligned}
\rt.\label{3.4}
\ee
and
\be
\lt\{
\begin{aligned}
&\wh^2_{tt}+|\xi|^2 \wh^2+\f{\mu}{(1+t)^\la}\wh^2_t=\hat{f},\\
&\wh(0,x)=0,\q \wh_t(0,x)=0.\label{3.5}
\end{aligned}
\rt.
\ee
The solution of \eqref{3.4} can be represented in the form
\be
\wh^1(t,\xi)=\Phi_1(t,\xi)\wh_0(\xi)+\Phi_2(t,\xi)\wh_1(\xi). \nn
\ee
And
The solution of \eqref{3.5} can be represented in the form
\be
\wh^2(t,\xi)=\int^t_0\Phi_2(t-\tau,\xi)\hat{f}(\xi)d\xi. \nn
\ee
$\Phi_1(t,\xi)$ and $\Phi_2(t,\xi)$ are the fundamental solutions of $\eqref{3.4}_1$ satisfying
\[
\Phi_1(0,\xi)=1,\q \p_t\Phi_1(0,\xi)=0,
\]
\[
\Phi_2(0,\xi)=0,\q \p_t\Phi_2(0,\xi)=1.
\]
If we set
\be
K_i(t,x)=\mathscr{F}^{-1}_{\xi\rightarrow x}[\Phi_i(t,\xi)] \qq i=1,2.\nn
\ee
Then the solution of \eqref{3.3} can be represented by
\bea
w(t,x)&=&(K_1(t,\cdot)\ast w_0(\cdot))(x)+(K_2(t,\cdot)\ast w_1(\cdot))(x)\nn\\
&&+\int^t_0(K_2(t-\tau,\cdot)\ast f(\cdot))(x)d\tau.\label{3.6}
\eea

The asymptotic behavior of $\Phi_i(t,\xi)(i=1,2)$ with respect to time has already been investigated carefully in Wirth \cite{Wirth}. Here we will give the most important Theorem there for our use. Before that ,we need to introduce some notations.

In the space $(t,\xi)\in(\bR^+\times\bR^n)$, denote $Z_1$, $Z_2$ and $Z_3$ as follows.
\be
Z_{1}=\Big\{(t,\xi):|\xi|\leq \f{1}{4}\f{\mu}{(1+t)^\la}\Big\},\nn
\ee
\be
Z_2=\Big\{(t,\xi):|\xi|\geq\f{1}{4}\f{\mu}{(1+t)^\la}\Big\}\cap\Big\{(t,\xi):|\xi|\leq 1\Big\},\nn
\ee
\be
Z_3=\Big\{(t,\xi):|\xi|\geq\f{1}{4}\f{\mu}{(1+t)^\la}\Big\}\cap\Big\{(t,\xi):|\xi|\geq 1\Big\}.\nn
\ee
\begin{lemma}[Theorem17 of \cite{Wirth}] There exists a constant $C_0$, depending on $\la,\mu$ such that $\Phi_i$ has the following estimates in different zones

When $(t,\xi)\in Z_1$,
\be
|\Phi_i(t,\xi)|\leq C_0e^{-C_0|\xi|^2(1+t)^{1-\la}}. \label{3.7}
\ee

When $(t,\xi)\in Z_2$,

\be
|\Phi_i(t,\xi)|\leq C_0e^{-C_0(1+t)^{1-\la}}e^{C_0(1-|\xi|^2)(1+t_\xi)^{1-\la}}, \label{3.8}
\ee
where $t_\xi$ is the upper boundary line of the zone $Z_1$, which means
\be
|\xi|=\f{1}{4}\f{\mu}{(1+t_\xi)^\la}.  \label{3.9}
\ee

When $(t,\xi)\in Z_3$,
\be
|\Phi_i(t,\xi)|\leq C_0e^{-C_0(1+t)^{1-\la}}.\label{3.10}
\ee

\end{lemma}
\pf This is a direct result of Theorem17 in Wirth\cite{Wirth}. We omit the detail. \ef

Next we give two integral properties of $\Phi_i$ in the phase space $\xi.$
\begin{lemma}
\be
\|\xi^\al\Phi_i(t,\xi)\|_{L^1(Z_j)}\leq C(1+t)^{-(1-\la)(n/2+|\al|/2)},\label{3.11}
\ee
\be
\|\xi^\al\Phi_i(t,\xi)\|_{L^2(Z_j)}\leq C(1+t)^{-(1-\la)(n/4+|\al|/2)},\label{3.12}
\ee
where $i,j=1,2$.
\end{lemma}
\pf In the zone $Z_1$, using \eqref{3.7}, we have
\bea
&&\|\xi^\al\Phi_i(t,\xi)\|_{L^1(Z_1)}\nn\\
&\leq&\int\limits_{\xi\in Z_1} |\xi|^\al|\Phi_i(t,\xi)|d\xi\nn\\
&\leq&C\int^\i_0|\xi|^{n-1+\al} e^{-C_0|\xi|^2(1+t)^{1-\la}}d|\xi|\nn\\
&\leq&C(1+t)^{-(1-\la)(n/2+|\al|/2)}\int^\i_0 s^{n-1+|\al|}e^{-C_0s^2}ds\nn\\
&\leq&C(1+t)^{-(1-\la)(n/2+|\al|/2)}.\nn
\eea
Integration in $Z_2$ will use the representation of $t_\xi$. From \eqref{3.7}, we have $1+t_\xi=\big(\f{4}{\mu}|\xi|\big)^{-\f{1}{\la}}$. Then using \eqref{3.8}, we have
\bea
&&\|\xi^\al\Phi_i(t,\xi)\|_{L^1(Z_2)}\nn\\
&\leq& C_0e^{-C_0(1+t)^{1-\la}}\int\limits_{\xi\in Z_2}|\xi|^\al e^{C_0(1-|\xi|^2)\big(\f{4}{\mu}|\xi|\big)^{-\f{1-\la}{\la}}}d\xi\nn\\
&\leq&Ce^{-C_0(1+t)^{1-\la}}\int^1_{\f{\mu}{4}(1+t)^{-\la}}|\xi|^{n-1+|\al|}e^{C_0\big(\f{4}{\mu}|\xi|\big)^{-\f{1-\la}{\la}}}
e^{-C|\xi|^{2-\f{1-\la}{\la}}}d|\xi|\nn
\eea
Noting that
\bea
&&|\xi|^{n-1+|\al|}e^{-C|\xi|^{2-\f{1-\la}{\la}}}\nn\\
&=&\Big(|\xi|^{2-\f{1-\la}{\la}}\Big)^{\f{n-1+|\al|}{2}}e^{-C|\xi|^{2-\f{1-\la}{\la}}}|\xi|^{\f{1-\la}{2\la}(n-1+|\al|)}\nn\\
&\leq&C|\xi|^{\f{1-\la}{2\la}(n-1+|\al|)}. \nn
\eea
Then we get
\bea
&&\|\xi^\al\Phi_i(t,\xi)\|_{L^1(Z_2)}\nn\\
&\leq&Ce^{-C_0(1+t)^{1-\la}}\int^1_{\f{\mu}{4}(1+t)^{-\la}}|\xi|^{\f{1-\la}{2\la}(n-1+|\al|)}e^{C_0\big(\f{4}{\mu}|\xi|\big)^{-\f{1-\la}{\la}}}d|\xi|\nn\\
&=&Ce^{-C_0(1+t)^{1-\la}}\int^{(1+t)^{1-\la}}_{\big(\f{4}{\mu}\big)^{-{\f{1-\la}{\la}}}}s^{-\f{1}{2}(n-1+|\al|)-\f{1}{1-\la}}e^{C_0s}ds
\qq\Big(s=\big(\f{4}{\mu}|\xi|\big)^{-\f{1-\la}{\la}}\Big)\nn\\
&\leq&Ce^{-C_0(1+t)^{1-\la}}\underbrace{\int^{(1+t)^{1-\la}}_{\big(\f{4}{\mu}\big)^{-{\f{1-\la}{\la}}}}s^{-\f{1}{2}(n+|\al|)}e^{C_0s}ds}_I.\label{3.13}
\eea
$I$ can be dealt with integration by parts
\bea
I&=&\f{1}{C_0}s^{-\f{1}{2}(n+|\al|)}e^{C_0s}\bigg |^{s=(1+t)^{1-\la}}_{s=\big(\f{4}{\mu}\big)^{-{\f{1-\la}{\la}}}}\nn\\
&&-\f{n+|\al|}{2C_0}\int^{(1+t)^{1-\la}}_{\big(\f{4}{\mu}\big)^{-{\f{1-\la}{\la}}}}s^{-\f{1}{2}(n+|\al|)-1}e^{C_0s}ds\nn\\
&\leq&C(1+t)^{-\f{1-\la}{2}(n+|\al|)}e^{C_0(1+t)^{1-\la}}.\nn
\eea
Inserting this into \eqref{3.13}, we obtain that
\be
\|\xi^\al\Phi_i(t,\xi)\|_{L^1(Z_2)}\leq C(1+t)^{-(1-\la)(n/2+|\al|/2)}.\nn
\ee
The proof of \eqref{3.12} will be essentially the same as \eqref{3.11}. We omit the detail.

Combining the estimates of $\Phi_i(t,\xi)$ in Lemma3.1 and Lemma3.2, we have the following estimates to the solution of \eqref{3.2}.
\begin{lemma}
Assume that for a function $g(x)\in L^1\cap H^{s+l}(\bR^n)$, where $s=[n/2]+1$ and $l\geq 0$. Then we have
\be
\|\p^k_x(K_i(t,\cdot)\ast g(\cdot))\|_\i\leq C(1+t)^{-(1-\la)\f{n+k}{2}}\|g\|_1+Ce^{-C_0(1+t)^{1-\la}}\|\p^{s+k}_x g\|,\label{3.14}
\ee
where $k\in\bN$ and $0\leq k\leq l$.
\be
\|\p^k_x(K_i(t,\cdot)\ast g(\cdot))\|\leq C(1+t)^{-(1-\la)(\f{n}{4}+\f{k}{2})}\|g\|_1+Ce^{-C_0(1+t)^{1-\la}}\|\p^k_x g\|,\label{3.15}
\ee
where $k\in\bN$ and $0\leq k\leq s+l$.
\end{lemma}
\pf Using Lemma3.1, Lemma3.2 and H\"{o}lder inequality, we have
\bea
&&\|\p^k_x(K_i(t,\cdot)\ast g(\cdot))\|_{\i}\nn\\
&\leq&\||\xi|^k\hat{K_i}(t,\xi)\hat{g}(\xi)\|_1\nn\\
&=&\Big(\int\limits_{\xi\in Z_1\cup Z_2}+\int\limits_{\xi\in Z_3}\Big)|\xi|^k|\Phi_i(t,\xi)||\hat{g}|d\xi\nn\\
&\leq& \|\hat{g}\|_\i\int\limits_{\xi\in Z_1\cup Z_2}|\xi|^k|\Phi_i(t,\xi)|d\xi+\int\limits_{\xi\in Z_3}\|\Phi_i(t,\xi)\|_\i|\xi|^k|\hat{g}|d\xi\nn\\
&\leq& \|\hat{g}\|_\i\int\limits_{\xi\in Z_1\cup Z_2}|\xi|^k|\Phi_i(t,\xi)|d\xi\nn\\
&&+Ce^{-C_0(1+t)^{1-\la}}\Big(\int_{|\xi|\geq 1}|\xi|^{-2([n/2]+1)}d\xi\Big)^{\f{1}{2}}\Big(\int_{|\xi|\geq 1}|\xi|^{2([n/2]+1+k)}|\hat{g}|^2d\xi\Big)^{\f{1}{2}}\nn\\
&\leq&C(1+t)^{-(1-\la)\f{n+k}{2}}\|g\|_1+Ce^{-C_0(1+t)^{1-\la}}\|\p^{k+s}_xg\|. \nn
\eea
This proves \eqref{3.14}.

Using Plancherel equality, Lemma3.1 and Lemma3.2, we have
\bea
&&\|\p^k_x(K_i(t,\cdot)\ast g(\cdot))\|^2\nn\\
&=&\||\xi|^k\hat{K_i}(t,\xi)\hat{g}(\xi)\|^2\nn\\
&=&\Big(\int\limits_{\xi\in Z_1\cup Z_2}+\int\limits_{\xi\in Z_3}\Big)|\xi|^{2k}|\Phi_i(t,\xi)|^2|\hat{g}|^2d\xi\nn\\
&\leq& \|\hat{g}\|^2_\i\int\limits_{\xi\in Z_1\cup Z_2}|\xi|^{2k}|\Phi_i(t,\xi)|^2d\xi+\int\limits_{\xi\in Z_3}\|\Phi_i(t,\xi)\|^2_\i|\xi|^{2k}|\hat{g}|^2d\xi\nn\\
&\leq& \|g\|^2_1\int\limits_{\xi\in Z_1\cup Z_2}|\xi|^{2k}|\Phi_i(t,\xi)|^2d\xi+Ce^{-2C_0(1+t)^{1-\la}}\int\limits_{\xi\in Z_3}|\xi|^{2k}|\hat{g}|^2d\xi\nn\\
&\leq&C(1+t)^{-(1-\la)(\f{n}{2}+k)}\|g\|^2_1+Ce^{-2C_0(1+t)^{1-\la}}\|\p^{k}_xg\|^2. \nn
\eea
This proves \eqref{3.15}.\ef
\subsection{Asymptotic behavior of $v$}

From \eqref{2.3} and \eqref{3.6}, we have
\bea
v(t,x)&=&(K_1(t,\cdot)\ast v_0(\cdot))(x)+(K_2(t,\cdot)\ast v_1(\cdot))(x)\nn\\
&&+\int^t_0(K_2(t-\tau,\cdot)\ast Q(v,u)(\tau,\cdot))(x)d\tau,\label{3.16}
\eea
where
\be
v_1(x)=\p_tv|_{t=0}=-\nabla\cdot u_0-u_0\cdot\nabla v_0-\f{\g-1}{2}v_0\nabla\cdot u_0. \label{3.17}
\ee

 Before we do estimates on $v$, we first estimate the nonlinear term $Q(v,u)$ under the help of Theorem1.1. Remember that
\bea
Q(v,u)&=&\f{\mu}{(1+t)^\la}(-u\cdot\nabla v-\f{\g-1}{2}v\nabla\cdot u)             \nn            \\
&&-\p_t(u\cdot\nabla v-\f{\g-1}{2}v\nabla\cdot u)+\nabla\cdot(u\cdot\nabla u+\f{\g-1}{2}v\nabla v).  \label{3.18}
\eea
\begin{proposition}
Assume $f,g\in H^l(\bR^n)$. Then for $|\al|\leq l$, we have
\be
\|\p^\al_x(fg)\|\leq C_l\big(\|f\|_\i\|\p^l_xg\|+\|g\|_\i\|\p^l_xf\|\big).\nn
\ee
\end{proposition}

Proof of the proposition can be found in many literatures. See \cite{KM} for example.
\begin{lemma}
Under the assumption of Theorem1.1, we have the following estimates for $Q$.
\be
\|Q(v,u)(t,\cdot)\|_1\leq C\ve^2(1+t)^{-B-\f{1+\la}{2}}. \label{3.19}
\ee
And
\be
\|\p^k_xQ(v,u)(t,\cdot)\|\leq C\ve^2(1+t)^{-B-\f{1+\la}{2}}, \label{3.20}
\ee
where $k=0,1,2,...,s+m-2$.
\end{lemma}
\pf This is a direct computation of the combination of Leibniz formula, H\"{o}lder inequality, Proposition3.1 and \eqref{1.4} in Theorem1.1.\ef
\begin{lemma}
Suppose that $a>1$ and $a\geq b>0$. Then there exists a constant $C$ such that for all $t\geq0$,
\be
\int^t_0(1+t-\tau)^{-a}(1+\tau)^{-b}d\tau\leq C(1+t)^{-b}. \label{3.21}
\ee
\end{lemma}
\pf This is just a direct computation
\bea
&&\int^t_0(1+t-\tau)^{-a}(1+\tau)^{-b}d\tau\nn\\
&=&\Big(\int^{\f{t}{2}}_0+\int^t_{\f{t}{2}}\Big)(1+t-\tau)^{-a}(1+\tau)^{-b}d\tau\nn\\
&\leq&C(1+t)^{-a}\int^{\f{t}{2}}_0(1+\tau)^{-b}d\tau+C(1+t)^{-b}\int^t_{\f{t}{2}}(1+t-\tau)^{-a}d\tau\nn\\
&\leq&C(1+t)^{-a+1-b}+C(1+t)^{-b}\nn\\
&\leq&C(1+t)^{-b}.\nn
\eea
\ef
\textbf{$L^\i$ estimate}

From \eqref{3.16}, using \eqref{3.14}, for $k=0,1,...,m-2$
\bea
&&\|\p^k_xv\|_\i\nn\\
&\leq&C(1+t)^{-(1-\la)\f{n+k}{2}}\|(v_0,v_1)\|_1+Ce^{-C_0(1+t)^{1-\la}}\|\p^{s+k}_x (v_0,v_1)\|\nn\\
&&+C\int^t_0(1+t-\tau)^{-(1-\la)\f{n+k}{2}}\|Q(v,u)\|_1d\tau\nn\\
&&+C\int^t_0e^{-C_0(1+t-\tau)^{1-\la}}\|\p^{s+k}_x Q(v,u)\|d\tau\nn\\
&\leq&C(1+t)^{-(1-\la)\f{n+k}{2}}\|(v_0,u_0)\|_{H^{s+m-1}}\nn\\
&&+C\ve^2\int^t_0(1+t-\tau)^{-(1-\la)\f{n+k}{2}}(1+\tau)^{-B-\f{1+\la}{2}}d\tau.  \nn
\eea
 Noting $B=\f{(1+\la)n}{2}-\dl$, when $(1-\la)\f{n+k}{2}= B+\f{1+\la}{2}$, we have
 \be
 k=k_c=\f{1+\la}{1-\la}(n+1)-n-\f{2\dl}{1-\la}.\nn
 \ee
\textbf{Case 1:} $0\leq k\leq k_c$\\
Using Lemma3.5 and the initial data assumption, we have
\bea
&&\|\p^k_xv\|_\i\nn\\
&\leq&C\ve(1+t)^{-(1-\la)\f{n+k}{2}}+C\ve^2(1+t)^{-(1-\la)\f{n+k}{2}}\nn\\
&\leq&C\ve(1+t)^{-(1-\la)\f{n+k}{2}}. \label{3.22}
\eea
\textbf{Case 2:} $k_c\leq k \leq m-2$\\
\bea
&&\|\p^k_xv\|_\i\nn\\
&\leq&C\ve(1+t)^{-(1-\la)\f{n+k}{2}}+C\ve^2(1+t)^{-B-\f{1+\la}{2}}\nn\\
&\leq&C(1+t)^{-(1+\la)\f{n+1}{2}+\dl}. \label{3.23}
\eea
\textbf{$L^2$ estimate}

From \eqref{3.16}, using \eqref{3.15}, for $k=0,1,...,s+m-2$
\bea
&&\|\p^k_xv\|\nn\\
&\leq&C(1+t)^{-(1-\la)(\f{n}{4}+\f{k}{2})}\|(v_0,v_1)\|_1+Ce^{-C_0(1+t)^{1-\la}}\|\p^{k}_x (v_0,v_1)\|\nn\\
&&+C\int^t_0(1+t-\tau)^{-(1-\la)(\f{n}{4}+\f{k}{2})}\|Q(v,u)\|_1d\tau\nn\\
&&+C\int^t_0e^{-C_0(1+t-\tau)^{1-\la}}\|\p^{k}_x Q(v,u)\|d\tau\nn\\
&\leq&C(1+t)^{-(1-\la)(\f{n}{4}+\f{k}{2})}\|(v_0,u_0)\|_{H^{s+m-1}}\nn\\
&&+C\ve^2\int^t_0(1+t-\tau)^{-(1-\la)(\f{n}{4}+\f{k}{2})}(1+\tau)^{-B-\f{1+\la}{2}}d\tau.  \nn
\eea
When $(1-\la)(\f{n}{4}+\f{k}{2})=B+\f{1+\la}{2}$, we have
\be
k=k_c+\f{n}{2}.\nn
\ee
\textbf{Case 1:} $0\leq k\leq k_c+\f{n}{2}$\\
Using Lemma3.5, we have
\bea
&&\|\p^k_xv\|\nn\\
&\leq&C\ve(1+t)^{-(1-\la)(\f{n}{4}+\f{k}{2})}+C\ve^2(1+t)^{-(1-\la)(\f{n}{4}+\f{k}{2})}\nn\\
&\leq&C\ve(1+t)^{-(1-\la)(\f{n}{4}+\f{k}{2})}. \label{3.24}
\eea
\textbf{Case 2:} $k_c+\f{n}{2}\leq k\leq m-2$\\
\bea
&&\|\p^k_xv\|\nn\\
&\leq&C\ve(1+t)^{-(1-\la)(\f{n}{4}+\f{k}{2})}+C\ve^2(1+t)^{-B-\f{1+\la}{2}}\nn\\
&\leq&C\ve(1+t)^{-(1+\la)\f{n+1}{2}+\dl}.\label{3.25}
\eea
Combining the estimates \eqref{3.22}-\eqref{3.25}, we proved \eqref{1.5} and \eqref{1.6}.
\subsection{Asymptotic behavior of u}
Denote $u=(u^1,...,u^n)$. From $\eqref{2.2}_2$, differentiating it $k$ time in $x$, we have
\be
\p_t \p^k_xu^i+\f{\mu}{(1+t)^\la}\p^k_xu^i=-\p^k_x\big(\p_{i}v+u^j\p_{j}u^i+\f{\g-1}{2}v\p_{i}v\big). \label{3.26}
\ee
\textbf{$L^\i$ estimate}
\\
From \eqref{3.26}, we have
\be
\f{d}{dt}\Big[e^{\f{\mu}{1-\la}(1+t)^{1-\la}}\p^k_xu^i\Big]=
-e^{\f{\mu}{1-\la}(1+t)^{1-\la}}\p^k_x\big(\p_{i}v+u^j\p_{j}u^i+\f{\g-1}{2}v\p_{i}v\big). \label{3.27}
\ee
Integrating \eqref{3.27} from 0 to $t$, we obtain
\bea
&&e^{\f{\mu}{1-\la}(1+t)^{1-\la}}\p^k_xu^i=e^{\f{\mu}{1-\la}}\p^k_xu^i_0(x)\nn\\
&&\qq\qq\qq-\int^t_0e^{\f{\mu}{1-\la}(1+\tau)^{1-\la}}\p^k_x\big(\p_{i}v+u^j\p_{j}u^i+\f{\g-1}{2}v\p_{i}v\big)d\tau. \label{3.28}
\eea
From the estimates of $v$ \eqref{1.5} and \eqref{1.6} and \eqref{1.4}, we have
\bea
&&|-\p^k_x\big(\p_{i}v+u^j\p_{x_j}u^i+\f{\g-1}{2}v\p_{i}v\big)|\nn\\
&\leq&\lt\{
\begin{aligned}
&(1+t)^{-(1-\la)\f{n+k+1}{2}}\qq 0\leq k\leq k_c-1,  \\
&(1+t)^{-\f{1+\la}{2}(n+1)+\dl}\qq k_c-1<k\leq m-3.
\end{aligned}
\rt.  \label{3.29}
\eea
Then from \eqref{3.28} and \eqref{3.29}, we have \\
\textbf{Case 1:} $0\leq k\leq k_c-1$
\bea
\|\p^k_xu^i\|_\i&\leq& C\ve e^{-\f{\mu}{1-\la}(1+t)^{1-\la}}\nn\\
&&+e^{-\f{\mu}{1-\la}(1+t)^{1-\la}}\underbrace{\int^t_0(1+\tau)^{-(1-\la)\f{n+k+1}{2}}e^{\f{\mu}{1-\la}(1+\tau)^{1-\la}}d\tau}_I. \label{3.30}
\eea
Using integration by parts, it is easy to see
\bea
I&=&C\int^t_0(1+\tau)^{-(1-\la)\f{n+k+1}{2}+\la}de^{\f{\mu}{1-\la}(1+\tau)^{1-\la}}\nn\\
&\leq& C(1+t)^{-(1-\la)\f{n+k+1}{2}+\la}e^{\f{\mu}{1-\la}(1+t)^{1-\la}}.\nn
\eea
So, from \eqref{3.30}, we get
\be
\|\p^k_xu^i\|_\i\leq C\ve(1+t)^{-(1-\la)\f{n+k+1}{2}+\la}. \label{3.31}
\ee
\textbf{Case 2:}  $k_c-1\leq k\leq m-3 $
\bea
\|\p^k_xu^i\|_\i&\leq& C\ve e^{-\f{\mu}{1-\la}(1+t)^{1-\la}}\nn\\
&&+e^{-\f{\mu}{1-\la}(1+t)^{1-\la}}\int^t_0(1+\tau)^{-(1+\la)\f{n+1}{2}+\dl}e^{\f{\mu}{1-\la}(1+\tau)^{1-\la}}d\tau\nn\\
&\leq&C\ve(1+t)^{-(1+\la)\f{n+1}{2}+\la+\dl}.\label{3.32}
\eea
\textbf{$L^2$ estimate}
\\
Multiplying \eqref{3.26} by $\p^k_xu^i$, then integrating it on $\bR^n$ and using H\"{o}lder inequality, we get
\bea
&&\f{d}{dt}\|\p^k_xu^i\|^2+\f{\mu}{(1+t)^\la}\|\p^k_xu^i\|^2\nn\\
&=&-\int_{\bR^n}\p^k_xu^i\p^k_x\big(\p_{i}v+u^j\p_{j}u^i+\f{\g-1}{2}v\p_{i}v\big)dx\nn\\
&\leq& \f{\mu}{2(1+t)^\la}\|\p^k_xu^i\|^2+C(1+t)^\la\|\p^k_x\p_{i}v\|^2\nn\\
&&+C(1+t)^\la\Big(\|\p^k_x(u^j\p_{j}u^i)\|^2+\|\p^k_x(v\p_{i}v)\|^2\Big)\nn\\
&\leq&\f{\mu}{2(1+t)^\la}\|\p^k_xu^i\|^2+C(1+t)^\la\|\p^k_x\p_{i}v\|^2\nn\\
&&+C\ve^4(1+t)^{-2B-1}. \label{3.33}
\eea
Form \eqref{1.5} and \eqref{1.6}, we have
\bea
&&\|\p^k_x\p_{i}v\|\nn\\
&\leq&\lt\{
\begin{aligned}
&C\ve(1+t)^{-(1-\la)(\f{n}{4}+\f{k+1}{2})}\qq 0\leq k\leq k_c+\f{n}{2}-1,\nn\\
&C\ve(1+t)^{-(1+\la)\f{n+1}{2}+\dl} \qq k_c+\f{n}{2}-1<k\leq s+m-3.
\end{aligned}
\rt.
\eea
\textbf{Case 1:} $0\leq k\leq k_c+\f{n}{2}-1$
\bea
&&\f{d}{dt}\|\p^k_xu^i\|^2+\f{\mu}{2(1+t)^\la}\|\p^k_xu^i\|^2\nn\\
&\leq&C\ve^2(1+t)^{-(1-\la)(\f{n}{2}+k+1)+\la}.\nn
\eea
Then we have
\bea
&&\f{d}{dt}\Big(\|\p^k_xu^i\|^2e^{\f{\mu}{2(1-\la)}(1+t)^{1-\la}}\Big)\nn\\
&\leq&C\ve^2\int^t_0(1+\tau)^{-(1-\la)(\f{n}{2}+k+1)+\la}e^{\f{\mu}{2(1-\la)}(1+\tau)^{1-\la}}d\tau  \label{3.34}
\eea
Integrating \eqref{3.34} from $0$ to $t$, we obtain
\be
\|\p^k_xu^i\|\leq C\ve(1+t)^{-(1-\la)(\f{n}{4}+\f{k+1}{2})+\la}. \label{3.35}
\ee
\textbf{Case 2:} when $k_c+\f{n}{2}-1<k\leq s+m-3$
\bea
&&\f{d}{dt}\|\p^k_xu^i\|^2+\f{\mu}{2(1+t)^\la}\|\p^k_xu^i\|^2\nn\\
&\leq&C\ve^2(1+t)^{-(1+\la)(n+1)+\la+2\dl}.\nn
\eea
The same as \eqref{3.34}, we have
\be
\|\p^k_xu^i\|\leq C\ve(1+t)^{-(1+\la)\f{n+1}{2}+\la+\dl}. \label{3.36}
\ee
Combining \eqref{3.31}, \eqref{3.32}, \eqref{3.35} and \eqref{3.36}, we proved \eqref{1.7} and \eqref{1.8}.
\section{Proof of Theorem1.3}

\q In this Section, we derive that the smooth solution of the Euler equations have a polynomially decayed lower bound in time while in three dimensions, the vorticity will decay exponentially.\\
\\
\textbf{Lower-bound decay rate of $(\rho,u)$}.

Define
\bea
&&F(t)=\int_{\bR^n}x\cdot(\rho u)dx,\qq M(t)=\int_{\bR^n}(\rho-1)dx, \nn\\
&&B(t)=\{x||x|\leq R+t\}.\nn
\eea
Due to the finite propagation of the solution and the compact support of the initial data, we have
\be
\supp (\rho(t)-1)\subseteq B(t),\qq\supp u\subseteq B(t).\nn
\ee
From $\eqref{1.1}_1$, we have
\be
\f{d}{dt}M(t)=\int_{\bR^n}\rho_t dx=-\int_{\bR^n}\nabla\cdot(\rho u)=0,\nn
\ee
which means $M(t)=M(0)$.
So using Cauchy-Schwartz inequality, we have
\be
q_0=M(0)\leq\Big|\int_{B(t)}(\rho-1)dx\Big|\leq C\|(\rho-1)\|(R+t)^{\f{n}{2}}. \label{4.1}
\ee
Using $\eqref{1.1}_2$ and integration by parts, we get
\be
F'(t)+\f{\mu}{(1+t)^\la}F(t)=\int_{\bR^n}\big(\rho|u|^2+n(p(\rho)-p(1))\big)dx. \label{4.2}
\ee
While noting the convexity of $p(\rho)=\f{1}{\g}\rho^\g,\g>1$ , so we have
\be
\int_{\bR^n}(p(\rho)-p(1))dx\geq\int_{\bR^n}p'(1)(\rho-1)dx=M(t)=M(0). \label{4.3}
\ee
Combining \eqref{4.2} and \eqref{4.3}, we have
\be
F'(t)+\f{\mu}{(1+t)^\la}F(t)\geq nq_0. \label{4.4}
\ee
Multiplying \eqref{4.4} by $e^{\f{\mu}{1-\la}(1+t)^{1-\la}}$ and integrating it on $[0,t]$, we obtain
\bea
F(t)&\geq& e^{\f{\mu}{1-\la}}F(0)e^{-\f{\mu}{1-\la}(1+t)^{1-\la}}\nn\\
&&+nq_0e^{-\f{\mu}{1-\la}(1+t)^{1-\la}}\int^t_0e^{\f{\mu}{1-\la}(1+\tau)^{1-\la}}d\tau\nn\\
&\geq& e^{\f{\mu}{1-\la}}F(0)e^{-\f{\mu}{1-\la}(1+t)^{1-\la}}\nn\\
&&+\f{nq_0}{\mu}e^{-\f{\mu}{1-\la}(1+t)^{1-\la}}\int^t_0\f{\mu}{(1+t)^\la}e^{\f{\mu}{1-\la}(1+\tau)^{1-\la}}d\tau\nn\\
&=&e^{\f{\mu}{1-\la}}F(0)e^{-\f{\mu}{1-\la}(1+t)^{1-\la}}\nn\\
&&+\f{nq_0}{\mu}e^{-\f{\mu}{1-\la}(1+t)^{1-\la}}\Big(e^{\f{\mu}{1-\la}(1+t)^{1-\la}}-e^{\f{\mu}{1-\la}}\Big)\nn\\
&\geq&Cq_0.\nn
\eea
when $t>t_0$, where $t_0$ is a suitably large constant depending on $\la,\mu$.

Using Cauchy-Schwartz inequality, finite propagation speed, and $|\rho|$ is uniformly bounded, we have
\be
Cq_0\leq F(t)\leq C\Big(\int_{B(t)}\rho^2|x|^2dx\Big)^{\f{1}{2}}\Big(\int|u|^2dx\Big)^{\f{1}{2}}\leq C(R+t)^{\f{n+2}{2}}\|u\|.\label{4.5}
\ee
\eqref{4.1} and \eqref{4.5} imply \eqref{1.9}.\\
\\
\textbf{Exponential decay of the vorticity.}

In three dimensions, the vorticity $\o$ of the equations satisfies
\be
\p_t \o+\f{\mu}{(1+t)^\la}\o+u\cdot\nabla \o-\o\nabla u=0.  \label{4.6}
\ee
Multiplying \eqref{4.6} by $\o$ and integrating it on $\bR^n$, we get
\bea
\f{1}{2}\f{d}{dt}\int|w|^2dx+\f{\mu}{(1+t)^\la}\int|\o|^2dx&\leq& C\int(|\o|^2|\nabla u|+|\o\cdot\nabla u\cdot \o|dx)\nn\\
&\leq&C\|\nabla u\|_\i\int |\o|^2dx. \nn
\eea
From Theorem 1.1, we have $\|\nabla u\|_\i\leq C\ve(1+t)^{\f{B+1+\la}{2}}$. Noting $\f{B+1+\la}{2}\geq \la$ and $\ve$ is small, so we obtain
\be
\f{1}{2}\f{d}{dt}\int|w|^2dx+\f{\mu}{2(1+t)^\la}\int|\o|^2dx\leq 0.\nn
\ee
This implies the exponential decay of $\o$ \eqref{1.10}.

\section{Appendix}
\textbf{Proof of Lemma2.1}\\
Remember
\be
\psi(t,x)=a\f{|x|^2}{(1+t)^{1+\la}}, \qq a=\f{(1+\la)\mu}{8}\Big(1-\f{\dl}{(1+\la)n}\Big).\nn
\ee
Then
 \be
 \begin{aligned}
 &\psi_t=-(1+\la)\f{a|x|^2}{(1+t)^{2+\la}}=-\f{1+\la}{1+t}\psi, \\
 &\nabla\psi=\f{2ax}{(1+t)^{1+\la}},\q \Dl\psi=\f{2an}{(1+t)^{1+\la}}.
 \end{aligned} \label{5.1}
 \ee
And
\be
\f{|\nabla\psi|^2}{-\psi_t}=\f{4a}{1+\la}\f{1}{(1+t)^\la}=\f{1}{2}\Big(1-\f{\dl}{(1+\la)n}\Big)\f{\mu}{(1+t)^\la}, \label{5.2}
\ee
\be
\Dl\psi=\Big[\f{(1+\la)n}{4}-\f{\dl}{4}\Big]\f{\mu}{(1+t)^{1+\la}}.  \label{5.3}
\ee
 Multiplying \eqref{2.3} by $e^{2\psi}\p_tv$ and $e^{2\psi}v$, we have
\bea
&&\p_t\Big[\f{e^{2\psi}}{2}\big((\p_t v)^2+|\nabla v|^2\big)\Big]-\nabla\cdot(e^{2\psi}\p_t v\nabla v) \nn  \\
&&\q +\ew \Big(\f{\mu}{(1+t)^{\la}}-\f{|\nabla\psi|^2}{-\psi_t}-\psi_t\Big)(\p_t v)^2+\underbrace{\f{\ew}{-\psi_t}|\psi_t\nabla v-\p_t v\nabla\psi|^2}_{I_1}  \nn  \\
&=&\ew \p_t v Q(v,u),    \label{5.4}
\eea
and
\bea
&&\p_t\Big[e^{2\psi}\Big(v\p_t v+\f{\mu}{2(1+t)^\la}v^2\Big)\Big]-\nabla\cdot(e^{2\psi} v\nabla v) \nn  \\
&&\q +\ew \Big\{|\nabla v|^2+\Big(-\psi_t+\f{\la}{2(1+t)}\Big)\f{\mu}{(1+t)^\la}v^2+\underbrace{2v\nabla\psi\cdot\nabla v}_{I_2}-2\psi_tv\p_v-(\p_tv)^2\Big\} \nn  \\
&=&\ew v Q(v,u).    \label{5.5}
\eea
We estimate $I_1,I_2$ as follows:
\bea
I_1&\geq&\f{\ew}{-\psi_t}\Big((1-\dl_1)\psi^2_t|\nabla v|^2-(1/\dl_1-1)v^2_t|\nabla \psi|^2\Big)  \nn  \\
   &=&\ew\Big\{(1-\dl_1)(-\psi_t)|\nabla v|^2-\f{1}{2}\Big(1-\f{\dl}{(1+\la)n}\Big)(1/\dl_1-1)\f{\mu}{(1+t)^\la}v^2_t\Big\}.  \label{5.6}
\eea
Choosing $\dl_1$ close to 1 such that
\be
\f{1}{2}\Big(1-\f{\dl}{(1+\la)n}\Big)(1/\dl_1-1)\leq \f{\dl}{2(1+\la)n},  \label{5.7}
\ee
and
\bea
I_2&=&4\ew v\nabla v\cdot\nabla \psi-\ew\nabla v^2\cdot\nabla\psi  \nn  \\
&=& 4\ew v\nabla v\cdot\nabla \psi-\nabla\cdot(\ew v^2\nabla\psi)+2\ew v^2|\nabla\psi|^2+\ew(\Dl\psi)v^2.\label{5.8}
\eea
Then Inserting \eqref{5.6} and \eqref{5.8} into \eqref{5.4} and \eqref{5.5}, we have
\bea
&&\p_t\Big[\f{e^{2\psi}}{2}\big((\p_t v)^2+|\nabla v|^2\big)\Big]-\nabla\cdot(e^{2\psi}\p_t v\nabla v) \nn  \\
&&\q +\ew\Big\{ \Big(\f{\mu}{2(1+t)^{\la}}-\psi_t\Big)(v_t)^2+(1-\dl_1)(-\psi_t)|\nabla v|^2\Big\}  \nn  \\
&\leq&\ew \p_t v Q(v,u),    \label{5.9}
\eea
and
\bea
&&\p_t\Big[e^{2\psi}\Big(v\p_t v+\f{\mu}{2(1+t)^\la}v^2\Big)\Big]-\nabla\cdot\{(e^{2\psi}( v\nabla v+v^2\nabla\psi)\} \nn  \\
&&\q +\ew \Big\{\underbrace{|\nabla v|^2+4v\nabla\psi\cdot\nabla v+\Big(\f{\mu}{(1+t)^\la}(-\psi_t)+2|\nabla\psi|^2\Big)v^2}_{I_3}  \nn  \\
&&\qq+\Big(\la+\f{(1+\la)n}{2}-2\dl\Big)\f{\mu}{2(1+t)^{1+\la}}v^2-2\psi_tvv_t-v^2_t\Big\} \nn  \\
&=&\ew v Q(v,u).    \label{5.10}
\eea
By \eqref{5.2}, we have
\bea
I_3&=&|\nabla v|^2+4v\nabla\psi\cdot\nabla v+4\f{1-\dl/\big(2(1+\la)n\big)}{1-\dl/\big((1+\la)n\big)}v^2|\nabla \psi|^2  \nn  \\
   &\geq&(1-\dl_2)|\nabla v|^2+4\Big(\f{1-\dl/\big(2(1+\la)n\big)}{1-\dl/\big((1+\la)n\big)}-1/\dl_2\Big)|\nabla\psi|^2v^2  \nn
\eea
Choosing $\dl_2$ close to 1 can assure that for some $\dl_3,\dl_4>0$, we have
\bea
I_3 &\geq& \dl_3\Big(|\nabla v|^2+|\nabla\psi|^2v^2\Big)\nn\\
&\geq&\dl_4\Big(|\nabla v|^2+\f{\mu}{(1+t)^\la}(-\psi_t)v^2\Big)\label{5.11}
\eea
Inserting \eqref{5.11} into \eqref{5.10}, we get
\bea
&&\p_t\Big[e^{2\psi}\Big(v\p_t v+\f{\mu}{2(1+t)^\la}v^2\Big)\Big]-\nabla\cdot\{(e^{2\psi}( v\nabla v+v^2\nabla\psi)\} \nn  \\
&&\q +\ew\Big\{\dl_4\Big(|\nabla v|^2+\f{\mu}{(1+t)^\la}(-\psi_t)v^2\Big)  \nn  \\
&&\qq+\Big(\la+\f{(1+\la)n}{2}-2\dl\Big)\f{\mu}{2(1+t)^{1+\la}}v^2-2\psi_tvv_t-v^2_t\Big\} \nn  \\
&\leq&\ew v Q(v,u).    \label{5.12}
\eea
Integrating \eqref{5.12} on $\bR^n$, we have
\bea
&&\f{d}{dt}\int_{\bR^n}e^{2\psi}\Big(v\p_t v+\f{\mu}{2(1+t)^\la}v^2\Big)dx  \nn  \\
&&\q +\int_{\bR^n}\ew\Big\{\dl_4\Big(|\nabla v|^2+\f{\mu}{(1+t)^\la}(-\psi_t)v^2\Big)+\Big(\la+\f{(1+\la)n}{2}-2\dl\Big)\f{\mu}{2(1+t)^{1+\la}}v^2\nn  \\
&&\qq\qq -2\psi_tvv_t-v^2_t\Big\}dx  \nn  \\
&\leq& \int_{\bR^n}\ew v Q(v,u)dx.  \label{5.13}
\eea

To absorb the negative term $-v^2_t$, we integrate \eqref{5.9} in $\bR^n$ and multiply it by $(K+t)^\la$, where $K$ is a sufficiently large constant. Then we get
\bea
&&\f{d}{dt}\Big[(K+t)^\la\int_{\bR^n}\f{e^{2\psi}}{2}\big( v^2_t+|\nabla v|^2\big)dx\Big]-\f{\la}{(K+t)^{1-\la}}\int_{\bR^n}\f{e^{2\psi}}{2}\big( v^2_t+|\nabla v|^2\big)dx  \nn  \\
&&\q+\int_{\bR^n}\ew\Big\{ \Big(\f{\mu}{2}-\psi_t(K+t)^\la\Big)v^2_t+(1-\dl_1)(-\psi_t)(K+t)^\la|\nabla v|^2\Big\} \nn  \\
&\leq&(K+t)^\la\int_{\bR^n}\ew \p_t v Q(v,u)dx.   \label{5.14}
\eea
Now adding $\nu\cdot$\eqref{5.13} to \eqref{5.14},where $\nu$ is a sufficient small constant, we get
\bea
 &&\f{d}{dt}\int_{\bR^n}e^{2\psi}\Big\{\f{(K+t)^\la}{2}\big( v^2_t+|\nabla v|^2\big)+\underbrace{\nu vv_t}_{I_4}+\f{\nu\mu}{2(1+t)^\la}v^2\Big\}dx\nn  \\
&&\q+\int_{\bR^n}e^{2\psi}\Big\{\Big(\f{\mu}{2}-\nu+(-\psi_t)(K+t)^\la-\f{\la}{2(K+t)^{1-\la}}\Big)v^2_t\nn  \\
&&\qq\qq\qq +\Big((1-\dl_1)(-\psi_t)(K+t)^\la+\dl_4\nu-\f{\la}{2(K+t)^{1-\la}}\Big)|\nabla v|^2 \nn  \\
&&\qq\qq\qq +\dl_4\nu\f{\mu}{(1+t)^\la}(-\psi_t)v^2+\Big(\la+\f{(1+\la)n}{2}-\f{\dl}{2}\Big)\f{\nu\mu}{2(1+t)^{1+\la}}v^2 \nn  \\
&&\qq\qq\qq\underbrace{-2\nu\psi_tvv_t}_{I_5}\Big\}dx  \nn  \\
&\leq& (K+t)^\la\int_{\bR^n}\ew \p_t v Q(v,u)dx+\nu\int_{\bR^n}\ew v Q(v,u)dx.
\eea

The terms $I_4$ and $I_5$ can be absorbed by the other positive terms by applying the smallness of $\nu$ and largeness of $K$ and the following Cauchy-Schwartz inequality
\bea
|I_4|&\leq& \f{(K+t)^\la}{4}v^2_t+\f{\nu^2}{(K+t)^\la}v^2 \nn  \\
&\leq&\f{(K+t)^\la}{4}v^2_t+\f{\nu^2}{(1+t)^\la}v^2, \label{5.16}
\eea
and
\bea
|I_5|&\leq&\f{\nu\dl_4}{2}(-\psi_t)\f{\mu}{(1+t)^\la}v^2+\f{2\nu}{\mu\dl_4}(-\psi_t)(1+t)^\la v^2_t \nn  \\
&\leq&\f{\nu\dl_4}{2}(-\psi_t)\f{\mu}{(1+t)^\la}v^2+\f{2\nu}{\mu\dl_4}(-\psi_t)(K+t)^\la v^2_t. \label{5.17}
\eea
Denote
\be
E(t)=\int_{\bR^n}e^{2\psi}\Big\{\f{(K+t)^\la}{2}\big( v^2_t+|\nabla v|^2\big)+\underbrace{\nu vv_t}_{I_4}+\f{\nu\mu}{2(1+t)^\la}v^2\Big\}dx.\nn
\ee
When $\nu$ is sufficiently small and $K$ is large, using \eqref{5.16}, there exists a small constant $c_\dl$ such that
\bea
&&\f{d}{dt}E(t)+\int_{\bR^n}\Big(\la+\f{(1+\la)n}{2}-\f{\dl}{2}\Big)\f{\nu\mu}{2(1+t)^{1+\la}}v^2dx\nn\\
&&\qq\q\ +c_\dl\int_{\bR^n}e^{2\psi}\Big\{\Big(1+(-\psi_t)(K+t)^\la\Big)v^2_t\nn  \\
&&\qq\qq\qq\q\ +\Big(1+(-\psi_t)(K+t)^\la\Big)|\nabla v|^2\nn\\
&&\qq\qq\qq\q\ +(1+t)^{-\la}(-\psi_t)v^2\Big\}dxd\tau\nn\\
&\leq&(K+t)^\la\int_{\bR^n}\ew \p_t v Q(v,u)dx+\int_{\bR^n}\ew v Q(v,u)dx. \label{5.18}
\eea
Using \eqref{5.16}, we see that
\bea
&&\f{1}{4}(K+t)^\la \Big[J(t;v_t)+J(t;|\nabla v|)\Big]+(\f{\nu\mu}{2}-\nu^2){(1+t)^{-\la}}J(t;v) \nn  \\
&& \leq E(t)\leq\nn\\
&&\f{3}{4}(K+t)^\la \Big[J(t;v_t)+J(t;|\nabla v|)\Big]+(\f{\nu\mu}{2}+\nu^2)(1+t)^{-\la}J(t;v). \label{5.19}
\eea
Multiplying \eqref{5.18} by $(K+t)^{B+\la}$ and using \eqref{5.19}, we have
\bea
&&\f{d}{dt}\Big[(K+t)^{B+\la} E(t)\Big]\nn\\
&&+\Big(c_\dl-\f{3}{4}(B+\la)(K+t)^{\la-1}\Big)(K+t)^{B+\la} \Big[J(t;v_t)+J(t;|\nabla v|)\Big]\nn\\
&&+\f{(K+t)^{B+\la}}{2(1+t)^{1+\la}}\bigg(\Big(\la+\f{(1+\la)n}{2}-\f{\dl}{2}\Big)\nu\mu-2(B+\la)(\f{\nu\mu}{2}+\nu^2)\bigg)J(t;v) \nn  \\
&&+c_\dl(K+t)^{B+2\la}\Big[J_\psi(v_t)+J_\psi(t;|\nabla v|)\Big]+c_\dl(K+t)^BJ_\psi(t;v)  \nn  \\
&\leq& \int_{\bR^n}\ew \bigg[\Big((K+t)^{B+2\la}v_t+\nu (K+t)^{B+\la} v\Big) Q(v,u)\bigg]dx. \label{5.20}
\eea

Integrating \eqref{5.20} over $[0,t]$, we can show that by choosing small $\nu$ and large $K$,  there exists a constant $C_0$ depending on $\la,\mu, \dl, R$ such that
\bea
&&(K+t)^{B+2\la}\Big[J(t;v_t)+J(t;|\nabla v|)\Big]+(K+t)^BJ(t;v)  \nn  \\
&&+\int^t_0(K+\tau)^{B+\la} \Big[J(\tau;v_\tau)+J(\tau;|\nabla v|)\Big]d\tau  \nn  \\
&&+\int^t_0(K+\tau)^{B+2\la}\Big[J_\psi(\tau;v_\tau)+J_\psi(\tau;|\nabla v|)\Big]d\tau\nn\\
&&+\int^t_0\Big[(K+\tau)^{B}J_\psi(\tau;v)+(K+\tau)^{B-1}J(\tau;v)\Big]d\tau  \nn  \\
&\leq&C\|(v(0),v_t(0),\p_xv(0))\|\nn\\
&&+C\int^t_0\int_{\bR^n}\ew \Big\{\Big((K+\tau)^{B+2\la}v_\tau+ (K+\tau)^{B+\la} v\Big) Q(v,u)\Big\}dxd\tau. \label{5.21}
\eea

Considering that $\int^t_0(K+\tau)^{B+\la} \big[J(\tau;v_\tau)+J(\tau;|\nabla v|)\big]d\tau$ has been estimated, we multiply \eqref{5.14} by $(K+t)^{B+1}$ and integrate it over $[0,t]$ to obtain
\bea
&&(K+t)^{B+1+\la}\Big(J(t;v_t)+J(t;|\nabla v|)\Big)\nn\\
&&-\int^t_0(\la+B+1)(K+\tau)^{B+\la}\big[J(\tau;v_\tau)+J(\tau;|\nabla v|)\big]d\tau   \nn \\
&&+\int^t_0(K+\tau)^{B+1+\la}[J_\psi(\tau;v_\tau)+J_\psi(\tau;|\nabla v|)\big]d\tau\nn\\
&&+\int^t_0(K+\tau)^{B+1}J(\tau;v_\tau)d\tau  \nn  \\
&\leq&C\|(v_t(0),\p_xv(0))\|+ C\int^t_0\int_{\bR^n}(K+\tau)^{B+1+\la}\ew \p_\tau v Q(v,u)dxd\tau. \label{5.22}
\eea
For small $\nu$, adding $\nu\cdot$\eqref{5.22} to \eqref{5.21}, we have
\bea
&&(K+t)^{B+1+\la}\Big[J(t;v_\tau)+J(t;|\nabla v|)\Big]+(K+t)^BJ(t;v) \nn  \\
&&\qq+\int^t_0(K+\tau)^{B+1+\la}[J_\psi(\tau;v_\tau)+J_\psi(\tau;|\nabla v|)\big]d\tau\nn\\
&&\qq+\int^t_0\big[(K+\tau)^{B+1}J(\tau,v_\tau)+(K+\tau)^{B+\la} J(\tau,|\nabla v|)\big]d\tau  \nn  \\
&&\qq+\int^t_0\big[(K+\tau)^{B}J_\psi(\tau,v)+(K+\tau)^{B-1}J(\tau,v)\big]d\tau \nn  \\
&\leq&C\|(v(0),v_t(0),\p_xv(0))\|\nn\\
&&+C\int^t_0\int_{\bR^n}\ew \Big\{\Big((K+\tau)^{B+1+\la}v_\tau+ (K+\tau)^{B+\la} v\Big) Q(v,u)\Big\}dxd\tau.\nn
\eea
This proves Lemma2.1.\ef

\indent

{\bf Acknowledgement.} I want to express my gratitude to my advisor Professor Huicheng Yin in Nanjing Normal University for his guidance about this work. This work is proceeded when I am visiting Department of Mathematics in University of California, Riverside. So, I also want to express my thanks to my co-advisor, Professor Qi S. Zhang, in UCR for his encouragement.

\
\\

\end{document}